\pdfoutput=1
\documentclass{elsarticle}
\usepackage[utf8]{inputenc}

\usepackage{amsfonts}
\usepackage[english]{babel}
\usepackage{amsthm}
\usepackage{amsmath}
\usepackage{amssymb}
\usepackage{tikz} 
\usepackage[overload]{empheq}
\usepackage[numbers]{natbib}
\usepackage[ruled,vlined,linesnumbered]{algorithm2e}
\DeclareMathAlphabet{\pazocal}{OMS}{zplm}{m}{n}
\usepackage [autostyle, english = american]{csquotes} 
\MakeOuterQuote{"}
\usepackage{enumitem}
\usepackage{xargs}
\usepackage{subfigure}
\usepackage{booktabs}   
\usepackage[backref=page]{hyperref}
\hypersetup{					 
   colorlinks,
   linkcolor={red},
   citecolor={blue},
   urlcolor={blue!80!black}
}
\usepackage{cleveref} 

\newtheorem{theorem}{Theorem}
\newtheorem{lemma}[theorem]{Lemma}
\newdefinition{remark}{Remark}
\newtheorem{definition}{Definition}
\newtheorem{corollary}{Corollary}
\newdefinition{example}{Example}
\newtheorem{proposition}{Proposition}
\newdefinition{question}{Question}

\usetikzlibrary{positioning,arrows,petri,calc}
\tikzset{
place/.style={
circle,
thick,
minimum size=6mm,
draw
},
transitionV/.style={
rectangle,
thick,
fill=black,
minimum height=6mm,
inner xsep=1pt
}
}

\newcommand{\graph}{\mathcal{G}}

\newcommand{\R}{\mathbb{R}}
\newcommand{\nat}{\mathbb{N}}
\newcommand{\nato}{\mathbb{N}_0}
\newcommand{\Rmax}{{\R}_{\text{\normalfont max}}}
\newcommand{\maxplus}{\oplus}
\newcommand{\maxtimes}{\otimes}

\newcommand{\Rbar}{\bar{\R}_{\text{\normalfont max}}}
\newcommand{\emptystr}{\mathsf{e}}
\newcommand{\arcs}{E}
\newcommand{\nodes}{N}
\newcommand{\nonegset}{\Gamma}
\newcommand{\nonegsetm}{\Gamma_{M}}
\newcommand{\xy}{\pazocal{S}}
\newcommand{\xx}{\pazocal{P}}
\newcommand{\yy}{\pazocal{I}}
\newcommand{\ltm}{\mu}
\newcommand{\lang}{\pazocal{L}}

\newcommand{\MP}{P}
\newcommand{\MI}{I}
\newcommand{\MC}{C}
\newcommand{\MA}{A}
\newcommand{\MB}{B}
\newcommand{\graphPIC}{\graph(\lambda \MP \oplus \lambda^{-1} \MI \oplus \MC)}
\newcommand{\graphPICm}{\graph(\lambda \MP, \lambda^{-1} \MI, \MC)}
\newcommand{\wA}{\mathsf{\MakeLowercase{\MA}}}
\newcommand{\wB}{\mathsf{\MakeLowercase{\MB}}}
\newcommand{\wP}{\mathsf{\MakeLowercase{\MP}}}
\newcommand{\wI}{\mathsf{\MakeLowercase{\MI}}}
\newcommand{\wC}{\mathsf{\MakeLowercase{\MC}}}
\newcommand{\wZ}{\mathsf{z}}
\newcommandx{\solPIC}[3][1=\MP,2=\MI,3=\MC]{\Lambda\left( #1, #2, #3 \right)}
\catcode`@=11
\def\plslash{\ifx\@currsize\normalsize
{\mathchoice
{\,\mbox{\raisebox{0.2ex}{$\scriptstyle\circ$}\kern-1ex$\setminus$}}
{\,\mbox{\raisebox{0.2ex}{$\scriptstyle\circ$}\kern-1ex$\setminus$}}%
{\,\mbox{\raisebox{0.14ex}{$\scriptscriptstyle\circ$}\kern-0.8ex%
${\scriptstyle\setminus}$}}%
{\,\mbox{\raisebox{0.14ex}{$\scriptscriptstyle\circ$}\kern-0.8ex%
${\scriptstyle\setminus}$}}}%
\else\ifx\@currsize\large\,\mbox{\raisebox{0.2ex}{$\scriptstyle\circ$}\kern-1ex$\setminus$}
\else\ifx\@currsize\small\,\mbox{\raisebox{0.2ex}{$\scriptstyle\circ$}\kern-1ex$\setminus$}
\else\,\mbox{\raisebox{0.2ex}{$\scriptstyle\circ$}\kern-0.1ex$\setminus$}
\fi\fi\fi}

\def\prslash{\ifx\@currsize\normalsize
{\mathchoice
{\mbox{\raisebox{0.2ex}{$\scriptstyle\circ$}\kern-1ex$/$}}
{\mbox{\raisebox{0.2ex}{$\scriptstyle\circ$}\kern-1ex$/$}}%
{\mbox{\raisebox{0.14ex}{$\scriptscriptstyle\circ$}\kern-0.8ex%
${\scriptstyle/}$}}%
{\mbox{\raisebox{0.14ex}{$\scriptscriptstyle\circ$}\kern-0.8ex%
${\scriptstyle/}$}}}%
\else\ifx\@currsize\large\mbox{\raisebox{0.2ex}{$\scriptstyle\circ$}\kern-1ex$/$}
\else\ifx\@currsize\small\mbox{\raisebox{0.2ex}{$\scriptstyle\circ$}\kern-1ex$/$}
\else\mbox{\raisebox{0.2ex}{$\scriptstyle\circ$}\kern-1ex$/$} \fi\fi\fi}

\def\plslashblack{\ifx\@currsize\normalsize
{\mathchoice
{\,\mbox{\raisebox{0.2ex}{$\scriptstyle\bullet$}\kern-1ex$\setminus$}}
{\,\mbox{\raisebox{0.2ex}{$\scriptstyle\bullet$}\kern-1ex$\setminus$}}%
{\,\mbox{\raisebox{0.14ex}{$\scriptscriptstyle\bullet$}\kern-0.8ex%
${\scriptstyle\setminus}$}}%
{\,\mbox{\raisebox{0.14ex}{$\scriptscriptstyle\bullet$}\kern-0.8ex%
${\scriptstyle\setminus}$}}}%
\else\ifx\@currsize\large\,\mbox{\raisebox{0.2ex}{$\scriptstyle\bullet$}\kern-1ex$\setminus$}
\else\ifx\@currsize\small\,\mbox{\raisebox{0.2ex}{$\scriptstyle\bullet$}\kern-1ex$\setminus$}
\else\,\mbox{\raisebox{0.2ex}{$\scriptstyle\bullet$}\kern-0.1ex$\setminus$}
\fi\fi\fi}

\def\prslash{\ifx\@currsize\normalsize
{\mathchoice
{\mbox{\raisebox{0.2ex}{$\scriptstyle\bullet$}\kern-1ex$/$}}
{\mbox{\raisebox{0.2ex}{$\scriptstyle\bullet$}\kern-1ex$/$}}%
{\mbox{\raisebox{0.14ex}{$\scriptscriptstyle\bullet$}\kern-0.8ex%
${\scriptstyle/}$}}%
{\mbox{\raisebox{0.14ex}{$\scriptscriptstyle\bullet$}\kern-0.8ex%
${\scriptstyle/}$}}}%
\else\ifx\@currsize\large\mbox{\raisebox{0.2ex}{$\scriptstyle\bullet$}\kern-1ex$/$}
\else\ifx\@currsize\small\mbox{\raisebox{0.2ex}{$\scriptstyle\bullet$}\kern-1ex$/$}
\else\mbox{\raisebox{0.2ex}{$\scriptstyle\bullet$}\kern-1ex$/$} \fi\fi\fi}
\catcode`@=12

\newcommand{\zorc}{}

\renewcommand*{\backref}[1]{}
\renewcommand*{\backrefalt}[4]{%
    \ifcase #1 (Not cited.)%
    \or        (Cited on page~#2.)%
    \else      (Cited on pages~#2.)%
    \fi}

\begin{document}

\title{The non-positive circuit weight problem\\in parametric graphs:\\a solution based on dioid theory}

\author[1]{Davide Zorzenon\corref{cor1}} 
\ead{zorzenon@control.tu-berlin.de}

\author[2]{Jan Komenda}
\ead{komenda@ipm.cz} 

\author[1,3]{J\"{o}rg Raisch}
\ead{raisch@control.tu-berlin.de}

\cortext[cor1]{Corresponding author}
\address[1]{Technische Universit\"at Berlin, Control Systems Group, Einsteinufer 17, D-10587 Berlin, Germany}
\address[2]{Institute of Mathematics, Czech Academy of Sciences, \v{Z}itn\'{a} 25, 115 67 Prague, Czech Republic}
\address[3]{Science of Intelligence, Research Cluster of Excellence, Marchstr. 23, 10587 Berlin, Germany}

\begin{abstract}
Let us consider a parametric weighted directed graph in which every arc $(j,i)$ has weight of the form $w((j,i))=\max(\MP_{ij}+\lambda,\MI_{ij}-\lambda,\MC_{ij})$, where $\lambda$ is a real parameter and $\MP$, $\MI$ and $\MC$ are arbitrary square matrices with elements in $\R\cup\{-\infty\}$.
In this paper, we design an algorithm that solves the \textit{Non-positive Circuit weight Problem} (NCP) on this class of parametric graphs, which consists in finding all values of $\lambda$ such that the graph does not contain circuits with positive weight.
This problem, which generalizes other instances of the NCP previously investigated in the literature, has applications in the consistency analysis of a class of discrete-event systems called P-time event graphs.
The proposed algorithm is based on max-plus algebra and formal languages, and improves the worst-case complexity of other existing approaches, achieving strongly polynomial time complexity $\pazocal{O}(n^4)$ (where $n$ is the number of nodes in the graph). 
\end{abstract}
\begin{keyword}
Parametric graphs; Non-positive circuit weight; Max-plus algebra; Formal languages; Dioid theory; P-time event graphs
\end{keyword}

\maketitle

\section{Introduction}

In graph theory, a classical problem is to check whether, given a weighted directed graph, there exists a circuit with positive (or negative) weight.
One of the simplest and most famous algorithms that solves this problem is due to Bellman and Ford~\cite{cormen2009graphalgorithms}.
The algebraic equivalent of this problem is related to linear inequalities in the tropical (max-plus or min-plus) algebra: given a square matrix $\MA$ in the max-plus algebra, the precedence graph $\graph(\MA)$ does not contain circuits with positive weight if and only if the max-plus inequality $x\succeq \MA \otimes x$ admits a real solution $x$.
The solution set of this kind of inequalities is often called zone or weighted digraph polyhedron~\cite{mine2004weakly,joswig2016weighted}.
We emphasize that the problem can be equivalently stated in the max-plus or in the min-plus algebra.

In the present paper, we consider parametric weighted directed graphs, in which weights of the arcs are variable and depend on some parameters.
In this context, we refer to the problem of finding all the values of the parameters such that the graph does not include circuits with positive weight as the \textit{Non-positive Circuit weight Problem} (NCP).
In particular, we are interested in studying this problem on a subclass of parametric weighted directed graphs whose weights depend only on one parameter $\lambda\in \R$.

It is known from the seminal work~\cite{karp1978characterization} of Karp that, in a graph with constant arc weights, there are no circuits with positive weight if and only if the maximum circuit mean of the graph is non-positive.
Based on this result, the NCP has been solved in literature in the cases when the weights of the arcs depend proportionally ($\graph(\lambda \MA)$) or inversely ($\graph(\lambda^{-1} \MA)$) on $\lambda$, in the max-plus sense (see, \textit{e.g.}, Theorem 1.6.18 in~\cite{butkovivc2010max}).
In standard algebra, this corresponds to having arc weights of the form, respectively, $w((j,i)) = \MA_{ij} + \lambda$ and $w((j,i)) = \MA_{ij} - \lambda$, where $w((j,i))$ indicates the weight of arc $(j,i)$ and $\MA_{ij}\in\R$.
Karp and Orlin provided two algorithms, running respectively in $\pazocal{O}(n^3)$ and $\pazocal{O}(n m \log n)$, that solve the NCP on parametric graphs with $n$ nodes and $m$ arcs, whose weights are of type $w((j,i)) = \MA_{ij}  + \MB_{ij}\times \lambda$, with $\MA_{ij}\in\R$ and $\MB_{ij}\in\{0,+1\}$~\cite{karp1981parametric}; a faster solution of the same problem, with time complexity $\pazocal{O}(n m+n^2 \log n)$, was found in~\cite{young1991faster}.
Moreover, Levner and Kats solved the NCP when $w((j,i)) = \MA_{ij} + \MB_{ij}\times \lambda$, with $\MA_{ij}\in\R$ and $\MB_{ij}\in\{-1,0,+1\}$ in strongly polynomial time complexity $\pazocal{O}(m n^2)$~\cite{levner1998parametric}.

We extend these problems to a class of more general parametric precedence graphs of the form $\graphPIC$, where $\MP,\ \MI$ and $\MC$ are three arbitrary square matrices of the same dimension in the max-plus algebra.
We refer to the NCP for this class of parametric precedence graphs as \textit{Proportional-Inverse-Constant}-NCP (PIC-NCP).
In standard algebra, the weight of a generic arc $(j,i)$ in graph $\graphPIC$ can be expressed as $w((j,i)) = \max(\MP_{ij} + \lambda, \MI_{ij} - \lambda, \MC_{ij})$, where $\MP_{ij},\MI_{ij},\MC_{ij}\in\R\cup\{-\infty\}$. 
Since the class of parametric graphs studied in the present paper includes all the above-mentioned ones, we aim to extend previous results on the NCP.

The interest in this problem comes from a class of discrete-event dynamical systems called P-time event graphs (P-TEGs)~\cite{CALVEZ19971487}.
In~\cite{zorzenon2021periodic}, it has been shown that, for all $d\in\nat$, a P-TEG 
with initially $0$ or $1$ token per place,
characterized by four square matrices $A^0,\ A^1,\ B^0,\ B^1$, admits consistent $d$-periodic trajectories of period $\lambda$ if and only if the precedence graph $\graph(\lambda B^{1\sharp} \oplus \lambda^{-1} (A^{1}\oplus E_{\otimes})  \oplus (A^{0} \oplus B^{0\sharp}))$ does not contain circuits with positive weight.
The problem of finding the periods of all admissible 1-periodic trajectories has been studied in~\cite{becha2013modelling,declerck2017extremum,lee2014steady} but an explicit formula for the admissible periods has not yet been found.
Moreover, in~\cite{zorzenon2020bounded}, we proved that a P-TEG is \textit{boundedly consistent} (i.e., there exists a consistent trajectory for the P-TEG in which the delay of the $k$-th firing of every pair of transitions is bounded for all $k$) if and only if it admits a 1-periodic trajectory.
Since in most P-TEG applications bounded consistency is not only desirable, but necessary, this further motivates our study.

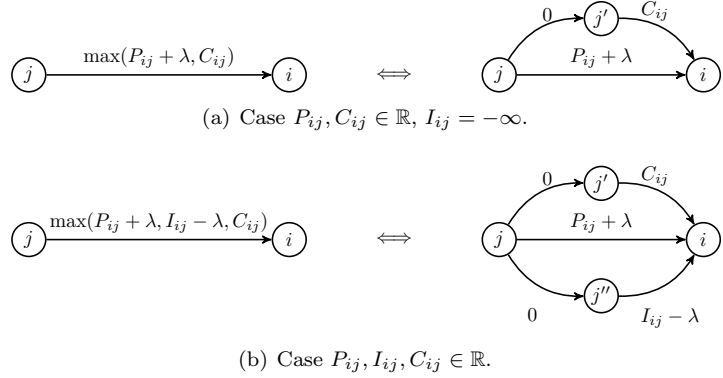
\begin{figure}[t]
\centering
\subfigure[Case $\MP_{ij},\MC_{ij}\in\R$, $\MI_{ij}=-\infty$.]{\label{fig:LevnerKatsa}
 \resizebox{.8\linewidth}{!}{
    \begin{tikzpicture}[node distance=3cm and 3cm,>=stealth',bend angle=45,thick]

\node [place,tokens=0,label=center:{$j$}] (P1) {};
\node [place,tokens=0,label=center:{$i$},right= 4 cm of P1] (P2) {};

\draw (P1) edge[->] node[auto,swap,label=above:{$\max(\MP_{ij} + \lambda,\MC_{ij})$}] {} (P2);
\node [right = 1cm of P2] (eq) {$\iff$};
\node [place,tokens = 0,label=center:{$j$},right = 1cm of eq] (P11) {};
\node [place,tokens=0,label=center:{$i$},right= of P11] (P22) {};
\node [place,tokens=0,label=center:{$j'$}] at ($(P11.east) + (1.5cm,1cm)$) (P33) {};
\draw (P11) edge[->] node[auto,swap,label=above:{$\MP_{ij} + \lambda$}] {} (P22);
\draw (P11) edge[bend left=30,->] node[auto,swap,label=above:{$0$}] {} (P33);
\draw (P33) edge[bend left=30,->] node[auto,swap,label=above:{$\MC_{ij}$}] {} (P22);


\end{tikzpicture}
    }}
\hfill
\subfigure[Case $\MP_{ij},\MI_{ij},\MC_{ij}\in\R$.]{\label{fig:LevnerKatsb}
  \resizebox{.8\linewidth}{!}{
    \begin{tikzpicture}[node distance=3cm and 3cm,>=stealth',bend angle=45,thick]

%

\node [place,tokens=0,label=center:{$j$}] (P1) {};
\node [place,tokens=0,label=center:{$i$},right= 4 cm of P1] (P2) {};
\draw (P1) edge[->] node[auto,swap,label=above:{$\max(\MP_{ij} + \lambda,\MI_{ij} - \lambda,\MC_{ij})$}] {} (P2);
\node [right = 1cm of P2] (eq) {$\iff$};
\node [place,tokens = 0,label=center:{$j$},right = 1cm of eq] (P11) {};
\node [place,tokens=0,label=center:{$i$},right= of P11] (P22) {};
\node [place,tokens=0,label=center:{$j'$}] at ($(P11.east) + (1.5cm,1cm)$) (P33) {};
\node [place,tokens=0,label=center:{$j''$}] at ($(P11.east) + (1.5cm,-1cm)$) (P44) {};
\draw (P11) edge[->] node[auto,swap,label=above:{$\MP_{ij} + \lambda$}] {} (P22);
\draw (P11) edge[bend left=30,->] node[auto,swap,label=above:{$0$}] {} (P33);
\draw (P33) edge[bend left=30,->] node[auto,swap,label=above:{$\MC_{ij}$}] {} (P22);
\draw (P11) edge[bend right=30,->] node[auto,swap,label=below:{$0$}] {} (P44);
\draw (P44) edge[bend right=30,->] node[auto,swap,label=below:{$\MI_{ij} - \lambda$}] {} (P22);

\end{tikzpicture}
    }}
\caption{Transformations needed to solve the PIC-NCP using Levner-Kats algorithm. Both transformations preserve the maximum weight of all paths from node $j$ to node $i$, for all values of $\MP_{ij},\MI_{ij},\MC_{ij},\lambda$.}
\label{fig:LevnerKats}
\end{figure}

We remark that the PIC-NCP can always be formulated as an instance of the NCP studied by Levner and Kats in~\cite{levner1998parametric}.
However, in order to do so, it is necessary to build an augmented precedence graph by adding $n'\in\{0,\ldots,2\times m\}$ nodes and $m'\in\{0,\ldots,4\times m\}$ arcs to $\graphPIC$, in particular: one node and two arcs for every arc $(j,i)$ in $\graphPIC$ for which two elements among $\MP_{ij}$, $\MI_{ij}$ and $\MC_{ij}$ are finite; two nodes and four arcs for every $(j,i)$ such that $\MP_{ij}$, $\MI_{ij}$ and $\MC_{ij}$ are all finite, as shown in Figure~\ref{fig:LevnerKats}.
This additional step increases the time complexity of the Levner-Kats algorithm to $\pazocal{O}((m+m')(n+n')^2)$, which leads to a worst-case complexity (attained when $n'=2\times m$ and $m=n^2$, i.e., when all elements of $\MP,\MI,\MC$ are real numbers) of $\pazocal{O}(n^6)$.

In this paper, we propose an algorithm based on techniques from dioid (or idempotent semiring) theory, which solves the PIC-NCP in time complexity $\pazocal{O}(n^4)$ and space complexity $\pazocal{O}(n^2)$, thus improving the worst-case complexity of the Levner-Kats algorithm.
Moreover, our algorithm provides a closed-formula expression for the extreme values of parameter $\lambda$ that solve the problem.
As a by-product, in this paper we indirectly prove that the class of linear programs of form~\eqref{lp:min} (or similarly~\eqref{lp:max}, see page~\pageref{lp:min}) can be solved in the same complexity.
The use of tropical algebra to solve linear programming problems is not new in literature, and is motivated by Smale's 9\textsuperscript{th} unsolved problem in mathematics~\cite{smale1998mathematical}, which asks whether linear programs admit a strongly polynomial time algorithm~\cite{butkovivc2019note,loho2016abstract,gaubert2010tropical}.

The algorithm is based on two different diods, the max-plus algebra and the semiring of formal languages.
The relation between matrices in the max-plus algebra and elements of a formal language is made explicit by means of \textit{multi--precedence graphs}, which are multi--directed graphs that generalize the concept of precedence graphs on multiple matrices.
Associating every matrix of a multi--precedence graph with a symbol, and every arc with an element of a matrix, we show how propositions on formal languages can be used to prove algebraic statements in the tropical semiring.

The paper is organized as follows.
In Section~\ref{se:alge}, some basic algebraic concepts from dioid theory (in particular, max-plus algebra and the algebra of formal languages) and weighted directed graph theory are recalled.
In Section~\ref{se:multiprecedencegraphs}, multi--precedence graphs are presented.
Section~\ref{se:main_problem} defines the main problem considered in the paper and gives its solution using both linear programming and a strongly polynomial algorithm based on dioid theory techniques, Algorithm~\ref{al:lambdasnn}.
In Section~\ref{se:validity}, the correctness of Algorithm~\ref{al:lambdasnn} is proven.
Finally, concluding remarks are stated in Section~\ref{se:conclusions}.

Notation: the set of positive, respectively non-negative, integers is denoted by $\nat$, respectively $\nato$.
Moreover, $\Rmax \coloneqq \R \cup \{-\infty\}$ and $\Rbar \coloneqq \Rmax \cup \{+\infty\}$.

\section{Preliminaries}\label{se:alge}

In this section, some basic concepts and results from dioid theory are summarized.
For more details, the reader is referred to~\cite{baccelli1992synchronization,heidergott2014max} and~\cite{hardouin2018control}.

\subsection{Dioid Theory}

A dioid $(\pazocal{D},\oplus,\otimes)$ is a set $\pazocal{D}$ endowed with two operations, $\oplus$ (addition) and $\otimes$ (multiplication), which have the following properties: both operations are associative and have a neutral element indicated, respectively, by $\epsilon$ (zero element) and $e$ (unit element), $\oplus$ is commutative and idempotent ($\forall a\in \pazocal{D}\ a\oplus a = a$), $\otimes$ distributes over $\oplus$, and $\epsilon$ is absorbing for $\otimes$ ($\forall a\in \pazocal{D}\ a\otimes \epsilon = \epsilon \otimes a = \epsilon$).
The operation $\oplus$ induces the natural order relation $\preceq$ on $\pazocal{D}$, defined by: $\forall a,b\in \pazocal{D}\ a \preceq b \Leftrightarrow a\oplus b = b$.
A dioid is complete if it is closed for infinite sums and if $\otimes$ distributes over infinite sums; in a complete dioid $(\pazocal{D},\oplus,\otimes)$, there exists a unique greatest (in the sense of $\preceq$) element of $\pazocal{D}$, denoted $\top$, \zorc{which is given by} $\top = \bigoplus_{x\in \pazocal{D}} x$.
The Kleene star of an element $a$ of a complete dioid, denoted $a^*$, is defined by $a^* = \bigoplus_{i\in \mathbb{N}_0} a^i$, with $a^0=e$ and $a^{i+1} = a\otimes a^i$. The operator $^+$ is defined as $a^+ = \bigoplus_{i\in \nat} a^i$; hence $a^* = e \oplus a^+$.
As in standard algebra, when unambiguous, the multiplication will be indicated simply as $a\otimes b = ab$.

If $(\pazocal{D},\oplus,\otimes)$ is a dioid, then the operations $\oplus$ and $\otimes$ can be extended to matrices with elements in $\pazocal{D}$: $\forall A,B\in \pazocal{D}^{m\times n},\ C\in\pazocal{D}^{n\times p}$
\[
	(A\oplus B)_{ij} = A_{ij}\oplus B_{ij},\quad  (A\otimes C)_{ij} = \bigoplus_{k=1}^n (A_{ik}\otimes C_{kj}).
\]
Moreover, the multiplication between a scalar and a matrix is defined as: $\forall  \lambda\in \pazocal{D},\ A\in\pazocal{D}^{m\times n}\quad (\lambda\otimes A)_{ij} = \lambda \otimes A_{ij}$.  
If $(\pazocal{D},\oplus,\otimes)$ is a complete dioid, then the set of $n\times n$ matrices endowed with $\oplus$ and $\otimes$ as defined above is a complete dioid, $(\pazocal{D}^{n\times n},\oplus,\otimes)$. Its zero and unit elements, respectively, are the matrices $\pazocal{E}$ and $E_\otimes$, where $\pazocal{E}_{ij} = \epsilon\ \forall i,j$ and $(E_\otimes)_{ij} = e$ if $i=j$, $(E_\otimes)_{ij} = \epsilon$ if $i\neq j$.
Furthermore, $A\preceq B \ \Leftrightarrow \ A_{ij} \preceq B_{ij} \ \forall i,j$.
We recall the following properties of the Kleene star operator $^*$.

\begin{proposition}[From~\cite{hardouin2018control}]
Let $\mathcal{D}$ be a complete dioid and $a,b\in\mathcal{D}$. The Kleene star operator $^*$ has the following properties:
\begin{gather}
\zorc{ a^*a^* = a^*} \label{eq:prop0} \\
(a\oplus b)^* = (a^* b)^* a^* = (b^* a)^* b^* \label{eq:prop1} \\
a (b a)^* = (a b)^* a \label{eq:prop2} \\
(a b^*)^* = e \oplus a (a\oplus b)^*. \label{eq:prop3}
\end{gather}
\end{proposition}

\subsection{Max-plus algebra}

An important example of a complete dioid is the max-plus algebra, $(\Rbar,\oplus,\linebreak \otimes)$, where $\oplus$ indicates the standard maximum operation, $\otimes$ indicates the standard addition, $\epsilon = -\infty$, $e=0$, $\top=+\infty$, and $\preceq$ coincides with the standard "less than or equal to".
For all $\lambda \in \R$, we indicate by $\lambda^{-1}$ the max-plus multiplicative inverse, i.e., $\lambda^{-1} \otimes \lambda = \lambda \otimes \lambda^{-1} = 0$. 
Note that, since $(\Rbar,\oplus,\otimes)$ is a complete dioid, the dioid $(\Rbar^{n\times n},\maxplus,\maxtimes)$ is also complete.
%

\subsection{Formal languages}

In this paper, it will be convenient to interpret (max,+) addition and multiplication of square matrices respectively as union and concatenation of formal languages.
This will be formally stated in Section~\ref{se:multiprecedencegraphs}.
A correspondence between (max,+) algebra and formal languages is justified by the fact that formal languages, endowed with the operations of union and concatenation, form a complete dioid~\cite{droste2009semirings}.

Let $\Sigma = \{\wA_1,\ldots, \wA_l\}$ be an \textit{alphabet} of $l$ \textit{letters}  (or \textit{symbols}) $\wA_1,\ldots,\wA_l$.
The set of all finite sequences of letters (or \textit{strings}) from $\Sigma$ is denoted by $\Sigma^*$.
A subset of $\Sigma^*$, $\lang \subseteq \Sigma^*$, is a \textit{formal language}, and its elements $s \in \lang$ are \textit{words}.
Let $\lang_1, \lang_2\subseteq \Sigma^*$, $\lang_1 = \{s_1,\ldots,s_n\},\ \lang_2 = \{t_1,\ldots,t_m\}$ be two languages; then $\lang_1+\lang_2 = \{s_1,\ldots,s_n,t_1,\ldots,t_m \}$ indicates the \textit{union} of the two languages, while $\lang_1\cdot \lang_2 =\lang_1 \lang_2= \{s_1t_1,s_1t_2,\ldots,s_1t_m,s_2t_1,s_2t_2,\ldots,s_2t_m,\ldots,s_nt_m\}$ indicates the language containing the \textit{concatenations} of all strings of $\lang_1$ and $\lang_2$.
Let $s \in \Sigma^*$, we will often indicate in the same way the (single string) language $s \coloneqq \{s\} \subseteq \Sigma^*$; the context will clarify whether we are referring to $s$ as a string or as a language.
Let $2^{\Sigma^*}$ indicate the power set of $\Sigma^*$ (the set of subsets of $\Sigma^*$).
Using the notation above, $(2^{\Sigma^*},+,\cdot )$ forms a complete dioid, with zero element the empty language $\emptyset=\{\}$ and unit element the language containing only the empty string $\emptystr$, i.e., $\{\emptystr\}$.
Note that, in contrast to standard notation, we denote the empty string by $\emptystr$.

We denote by $|s|$ the length of the word $s$ ($|\emptystr|=0$), and by $|s|_{\wA_i}$ the length of $s$ relative to letter $\wA_i$, i.e., the number of occurrences of the letter $\wA_i$ in $s$.
We indicate by $s(i)$, $1\leq i\leq |s|$, the $i$-th symbol of $s$.
\begin{definition}[Balanced string]
We say that a string $s$ is \textit{balanced} if $|s|_{\wA_1} = |s|_{\wA_2} = \ldots = |s|_{\wA_l}$.
Moreover, a string $s$ is \textit{$x*x$--balanced} if it is balanced and $s(1) = s(|s|)$, \textit{$x*y$--balanced} if it is balanced and $s(1) \neq s(|s|)$.
For convenience, we consider the empty string $\emptystr$ to be both $x*x$-- and $x*y$--balanced.
\end{definition}

\subsection{Precedence graphs}

\begin{definition}[Precedence graph]
	Let $\MA \in \Rmax^{n\times n}$.
	The \textit{precedence graph} associated with $\MA$ is the weighted directed graph $\mathcal{G}( \MA )=(\nodes,\arcs,w)$, where
	\begin{itemize}
		\item[-] $\nodes=\{1,\ldots,n\}$ is the set of nodes,
		\item[-] $\arcs\subseteq \nodes\times \nodes$ is the set of arcs, defined such that there is an arc  $(j,i)\in \arcs$ from node $j$ to node $i$ iff $ \MA _{ij}\neq -\infty$,
        \item[-] $w:\arcs\rightarrow \R$ is a function that associates a weight $w((j,i))= \MA _{ij}$ \zorc{with} every arc $(j,i)$ of $\mathcal{G}( \MA )$.
	\end{itemize}
	When matrix $ \MA $ depends on some real parameters, $ \MA  = \MA(\lambda_1,\ldots,\lambda_p)$, \linebreak$\lambda_1,\ldots,\lambda_p\in \R$, we say that $\graph( \MA )$ is a \textit{parametric precedence graph}.
\end{definition}
	A \textit{path} $\rho$ in $\mathcal{G}( \MA )=(\nodes,\arcs,w)$ is a sequence of nodes $(i_1,i_2,\ldots,i_{r+1})$, $r \geq 0$, \zorc{such that $(i_j,i_{j+1})\in \arcs$ for all $j=1,\ldots,r$ (i.e., with an arc from node $i_j$ to node $i_{j+1}$ for $j=1,\ldots,r$)}; 
we will use the notation
	\[
	\rho = i_1 \rightarrow i_2 \rightarrow \ldots \rightarrow i_{r+1}.
	\]
    The length of a path $\rho$, denoted by $|\rho|_L$, is the number of its arcs. Its weight, $|\rho|_W$, is the max-plus product (standard sum) of the weights of its arcs:
	\[
	|\rho|_L=r, \quad |\rho|_W= \bigotimes_{j=1}^{r}  \MA _{i_{j+1},i_j}.
	\]
	We define the weight of every path of length $|\rho|_L=0$ to be $|\rho|_W=0$.
	A path is \textit{elementary} if all its nodes are distinct.
	A path $\rho=i_1\rightarrow \ldots \rightarrow i_{r+1}$ is called \textit{circuit} if its initial and final nodes coincide, i.e., if $i_1=i_{r+1}$. A circuit $\rho=i_1\rightarrow \ldots\rightarrow i_{r+1}$ is called \textit{elementary} if the path $\tilde{\rho}=i_1\rightarrow \ldots\rightarrow i_{r}$ is elementary.

We recall from~\cite{heidergott2014max} that, given $\MA\in \Rmax^{n\times n}$, $r\in \nato$, $( \MA ^r)_{ij}$ is equal to the maximum weight of all paths in $\graph( \MA )$ from node $j$ to node $i$ of length $r$.
If there is no such path, then $( \MA^r)_{ij}=-\infty$.
We indicate by $\mbox{mcm}( \MA )$ the \textit{maximum circuit mean} of the precedence graph $\graph( \MA )$\footnote{\zorc{This quantity coincides with the largest max-plus eigenvalue of matrix $\MA$.}}, which can be computed as $\mbox{mcm}( \MA ) = \bigoplus_{k=1}^n (tr( \MA ^k))^{\frac{1}{k}}$, where $tr(M)=\bigoplus_{i=1}^{n} M_{ii}$ is the trace of matrix $M$, i.e., the max-plus sum of its diagonal elements, and $a^{\frac{1}{k}}$ represents, again in the max-plus sense, the $k$-th root of $a$, i.e., $(a^{\frac{1}{k}})^k=a$ \cite{baccelli1992synchronization}.

We indicate by $\nonegset$ the set of all precedence graphs that do not contain circuits with positive weight.
The following two propositions
connect max-plus algebra with precedence graphs.
    
\begin{proposition}[From~\cite{baccelli1992synchronization,gallai1958maximum} and Proposition 1.6.10 in~\cite{butkovivc2010max}] \label{pr:v_leq_Bv}
Let $\MA\in \Rmax^{n\times n}$.
Then, inequality $x \succeq  \MA  \otimes x$ has at least one solution $x\in \R^n$ if and only if $\mathcal{G}( \MA )$ does not contain any (elementary) circuit with positive weight, i.e., $\graph(\MA)\in\nonegset$.
\end{proposition}


\begin{proposition} \label{pr:aux}
Let $\MA \in \Rmax^{n\times n}$.
Then, the following statements are equivalent:
\begin{enumerate}[label={\normalfont(\roman*)}, ref=(\roman*)]
	\item\label{en:1} $\graph(\MA)\in\nonegset$,
	\item\label{en:2} for all $i\in\{1,\ldots,n\}$, $(\MA^*)_{ii} = 0$,
	\item\label{en:3} for all $i\in\{1,\ldots,n\}$, $(\MA^*)_{ii} \neq +\infty$.
\end{enumerate}
Moreover, there is a circuit with positive weight in $\graph(\MA)$ containing node $i$ if and only if $(\MA^*)_{ii} = +\infty$.

\end{proposition}
\begin{proof}
\ref{en:1} $\Leftrightarrow$ \ref{en:2}: the proof comes from the following equivalences: there are no circuits with positive weight in $\graph( \MA )$ $\Leftrightarrow$ $\forall i\in\{1,\ldots,n\},k\in\nato$ $ (\MA^k) _{ii}\preceq 0$ $\Leftrightarrow$ $ (\MA^*) _{ii} = 0 \oplus  \MA _{ii}\oplus  (\MA^2) _{ii} \oplus \ldots = 0$.

\ref{en:1} $\Rightarrow$ \ref{en:3}: obvious from \ref{en:1} $\Rightarrow$ \ref{en:2}.

\ref{en:1} $\Leftarrow$ \ref{en:3}: suppose that $\graph( \MA )\notin \nonegset$.
Then there is a circuit $\rho$ from some node $i$ of length $k$ such that $|\rho|_W = ( \MA ^k)_{ii}\succ0$.
Let us consider the circuit $\rho_h$, formed by going through $\rho$ $h\in\nat$ times.
We have that $|\rho_h|_L=k^h$ and $|\rho_h|_W = (|\rho|_W)^h\succ|\rho|_W$ (in standard notation, $|\rho_h|_W = h\times |\rho|_W$).
Note that, since $\MA^{k^h}_{ii}$ is equal to the greatest weight of all circuits including $i$ of length $k^h$, $( \MA ^{k^{h}})_{ii}  \succeq |\rho_h|_W$.
Therefore, $( \MA ^*)_{ii} = 0 \oplus  \MA _{ii} \oplus \ldots \oplus ( \MA ^k)_{ii} \oplus \ldots \oplus ( \MA ^{k^2})_{ii} \oplus \ldots \oplus ( \MA ^{k^3})_{ii}\oplus \ldots = +\infty$.

Regarding the last sentence of the proposition, the sufficiency comes from the proof of \ref{en:1} $\Leftarrow$ \ref{en:3}.
For the necessity, suppose that $(\MA^*)_{ii} = +\infty$ but there is no circuit with positive weight from node $i$.
Then, for all $k\in \nat$, $(\MA^k)_{ii} \preceq 0$.
Therefore, $(\MA^*)_{ii} = 0 \oplus \MA_{ii} \oplus (\MA^2)_{ii}\oplus \ldots = 0$, which contradicts the hypothesis.
\end{proof}

\zorc{
\begin{remark}
    From the latter proposition, we have that 
    \[
    \graph(\MA)\in\nonegset \ \Leftrightarrow \ tr(\MA ^*)=0 \qquad \mbox{and} \qquad \graph(\MA)\notin\nonegset \ \Leftrightarrow \ tr(\MA ^*)=+\infty \ .
    \]
\end{remark}
}



%

The problem of detecting the existence of circuits with positive weight in precedence graphs can be reduced to the problem of finding the shortest path from a node~\cite{cormen2009graphalgorithms}.
The shortest path from a node problem, as well as the computation of the maximum circuit mean, have been well studied in literature; the most classical algorithms for their solutions are, respectively, the Bellman-Ford algorithm and Karp's algorithm, both of strongly polynomial complexity $\mathcal{O}(n \times m)$ where $n$ is the number of nodes and $m$ is the number of edges~\cite{dasdan1998faster}.
Finally, we recall that, for all $\MA\in\Rbar^{n\times n}$ such that $\graph(\MA)\in\nonegset$, $\MA^*$ can be computed by using Floyd-Warshall algorithm of strongly polynomial complexity $\pazocal{O}(n^3)$~\cite{cormen2009graphalgorithms}.
In the case $\graph(\MA)\notin\nonegset$, the problem of computing $\MA^*$ is NP-hard~\cite{cormen2009graphalgorithms}.
Fortunately, the algorithm that will be presented in Section~\ref{se:main_problem} will never face this issue in practice; nevertheless, considering $\MA^*$ when $\graph(\MA)\notin\nonegset$ will be useful for some theoretical results.


\section{Multi--precedence graphs}\label{se:multiprecedencegraphs}

In this section, a new type of directed graph, called multi--precedence graph, is defined.
Multi--precedence graphs are a generalization of precedence graphs, in which every arc is labeled and each label corresponds to a different matrix.
Their definition is similar to the one of max-plus automata~\cite{gaubert1995performance}, with the difference that in multi--precedence graphs there are no initial and final states.
Moreover, max-plus automata are a modelling framework for dynamical systems, while multi--precedence graphs are used here as a tool to connect the concepts of precedence graphs and formal languages.

\begin{definition}[Multi--precedence graph]
Let $ \MA _1,\ldots, \MA _l$ be $n\times n$ matrices in $\Rmax$. The \textit{multi--precedence graph} associated with matrices $ \MA _1,\ldots, \MA _l$ is the weighted multi--directed graph $\mathcal{G}( \MA _1,\ldots, \MA _l)=(\nodes,\Sigma, \ltm ,\arcs)$, where
\begin{itemize}
	\item[-] $\nodes=\{1,\ldots,n\}$ is the set of nodes,
	\item[-] $\Sigma=\{\wA_1,\ldots,\wA_l\}$ is the alphabet of symbols $\wA_1,\ldots,\wA_l$,
    \item[-] $ \ltm :\Sigma\rightarrow \Rmax^{n\times n}$ is the morphism defined as
            \[
             \ltm (\wZ) = \begin{cases}
                                 \MA _1 & \mbox{if } \wZ=\wA_1 \\
                                \vdots & \\
                                 \MA _l & \mbox{if } \wZ=\wA_l,
                            \end{cases}
            \]
    \item[-] $\arcs\subseteq \nodes \times \nodes \times \Sigma$ is the set of labeled arcs, defined such that there is an arc $(j,i,\wZ)\in \arcs$ from node $j$ to node $i$ labeled $\wZ$ with weight $( \ltm (\wZ))_{ij}$ iff $( \ltm (\wZ))_{ij}\neq -\infty$.
\end{itemize}
When matrices $ \MA _1,\ldots, \MA _l$ depend on some real parameters, $\MA_1 = \MA_1 (\lambda_1,\dots,\lambda_p),\linebreak \ldots,\MA_l = \MA_l (\lambda_1,\ldots,\lambda_p)$, $\lambda_1,\ldots,\lambda_p\in \R$, we say that $\mathcal{G}( \MA _1,\ldots, \MA _l)$ is a \textit{parametric multi--precedence graph}.
\end{definition}

A \textit{path} in $\mathcal{G}( \MA _1,\ldots, \MA _l)$ is a sequence of alternating nodes and labels of the form $\sigma = (i_1,\wZ_1,i_2,\wZ_2,\ldots,\wZ_{r},i_{r+1})$, $r \geq 0$, \zorc{such that $(i_j,i_{j+1},\wZ_j)\in \arcs$ for all $j=1,\ldots,r$ (i.e., with an arc from node $i_j$ to node $i_{j+1}$ labeled $\wZ_j$ for $j=1,\ldots,r$)}; 
we will use the notation
\[
\sigma = i_1 \xrightarrow{\wZ_1} i_2 \xrightarrow{\wZ_2} \ldots \xrightarrow{\wZ_{r}} i_{r+1},
\]
and we will say that the path $\sigma$ is labeled $\wZ_r \wZ_{r-1} \cdots \wZ_1$.
The definitions of length of a path $|\sigma|_L$, weight of a path $|\sigma|_W$, elementary path and circuit are the same as for precedence graphs.
When not otherwise stated, we will use the convention to indicate matrices with uppercase letters ($ \MA , \MB ,C,\ldots$) and their associated symbols in a multi--precedence graph with the corresponding lowercase letters ($\graph( \MA , \MB ,C,\ldots) = (\nodes,\Sigma, \ltm ,\arcs)$ with $\Sigma = \{\wA, \wB,\mathsf{c},\ldots\}$ such that $ \ltm (\wA) =  \MA $, $ \ltm (\wB)= \MB $, $ \ltm (\mathsf{c}) = C,\ldots$).

\begin{example} \label{ex:1}

\begin{figure}[t]
\centering
\resizebox{.8\linewidth}{!}{
\begin{tikzpicture}[node distance=.5cm and 1cm,>=stealth',bend angle=45,thick]
\tikzset{
    between/.style args={#1 and #2}{
         at = ($(#1)!0.5!(#2)$)
    }
}

\node [place,tokens=0,label=center:{$1$}] (P1) {};
\node [place,tokens=0,label=center:{$2$},right= of P1] (P2) {};
\node [place,tokens=0,label=center:{$3$},right= of P2] (P3) {};
\node [left= 10pt of P1] (GP) {$\mathcal{G}(\MA)$}; 

\node [place,tokens=0,label=center:{$1$},right= 70pt of P3] (Q1) {};
\node [place,tokens=0,label=center:{$2$},right= of Q1] (Q2) {};
\node [place,tokens=0,label=center:{$3$},right= of Q2] (Q3) {};
\node [left= 10pt of Q1] (GQ) {$\mathcal{G}(\MB)$}; 

\node [below=70pt of GP] (point1) {};
\node [below=70pt of Q3] (point2) {};
\node [place,between=point1 and point2,tokens=0,label=center:{$2$},] (PQ2) {};
\node [place,tokens=0,label=center:{$1$},left= of PQ2] (PQ1) {};

\node [place,tokens=0,label=center:{$3$},right= of PQ2] (PQ3) {};
\node [left= 10pt of PQ1] (GPQ) {$\mathcal{G}(\MA,\MB)$}; 

%

\draw (P1) edge[->] node[auto,swap,label=above:{$2$}] {} (P2);
\draw (P3) edge[loop above,->] node {$-1$} (P3);

\draw (Q1) edge[bend left=60,->] node[auto,swap,label=above:{$0$}] {} (Q2);
\draw (Q2) edge[bend left=60,->] node[auto,swap,label=below:{$-2$}] {} (Q1);
\draw (Q2) edge[->] node[auto,swap,label=above:{$-2$}] {} (Q3);
\draw (Q3) edge[bend left=90,->] node[auto,swap,label=below:{$1$}] {} (Q1);

\draw (PQ1) edge[->] node[auto,swap,label=above:{$2,\wA$}] {} (PQ2);
\draw (PQ3) edge[loop above,->] node {$-1,\wA$} (PQ3);
\draw (PQ1) edge[bend left=60,->] node[auto,swap,label=above:{$0,\wB$}] {} (PQ2);
\draw (PQ2) edge[bend left=60,->] node[auto,swap,label=below:{$-2,\wB$}] {} (PQ1);
\draw (PQ2) edge[->] node[auto,swap,label=above:{$-2,\wB$}] {} (PQ3);
\draw (PQ3) edge[bend left=90,->] node[auto,swap,label=below:{$1,\wB$}] {} (PQ1);

\end{tikzpicture}
}
\caption{Precedence graphs $\mathcal{G}( \MA )$, $\mathcal{G}( \MB )$ and multi--precedence graph $\mathcal{G}( \MA , \MB )$ associated with matrices $ \MA $ and $ \MB $ of Example~\ref{ex:1}.}
\label{fig:matrixMultiplicationIntuition}
\end{figure}
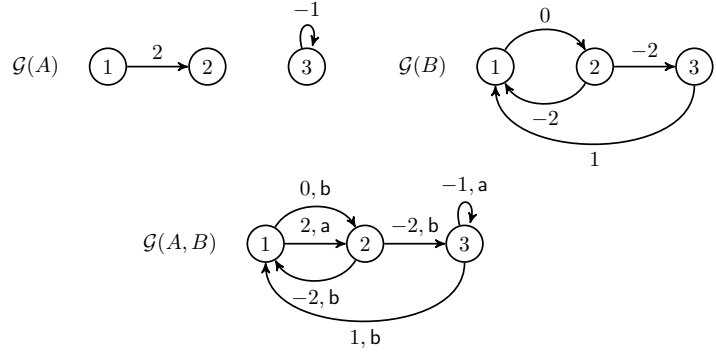

Let
\begin{equation*}
 \MA  =     
\begin{bmatrix}
-\infty & -\infty & -\infty \\
2 & -\infty & -\infty \\
-\infty & -\infty & -1
\end{bmatrix},~
 \MB  = \begin{bmatrix}
-\infty & -2 & 1 \\
0 & -\infty & -\infty \\
-\infty & -2 & -\infty
\end{bmatrix}.
\end{equation*}
The precedence graphs $\mathcal{G}( \MA )$ and $\mathcal{G}( \MB )$ and the multi--precedence graph $\mathcal{G}( \MA , \MB )$ are shown in Figure~\ref{fig:matrixMultiplicationIntuition}.
\end{example}

We can extend the morphism $ \ltm $ to $ \ltm :2^{\Sigma^*} \rightarrow \Rbar^{n\times n}$ as:
\[
 \ltm (\{\emptystr\}) = E_{\otimes}, \quad \ltm (\{ \wZ \}) = \mu( \wZ ),
\]
\[
 \ltm (\lang_1+\lang_2) =  \ltm (\lang_1) \oplus  \ltm (\lang_2),
\]
\[
 \ltm (\lang_1\lang_2) =  \ltm (\lang_1) \otimes  \ltm (\lang_2)
\]
for all $ \wZ \in \Sigma$, $\lang_1,\lang_2\subseteq \Sigma^*$.
It is trivial to see that the following properties hold for $ \ltm $.
Let $\lang_1,\lang_2\subseteq \Sigma^*$ be two languages such that $\lang_1\subseteq \lang_2$, then $ \ltm (\lang_1) \preceq  \ltm (\lang_2)$.
Moreover, for any language $\lang\subseteq \Sigma^*$, $\ltm(\lang^*) = \ltm(\lang)^*$.
If $s\in \Sigma^*$ is a string, we will often use the notation $\ltm(s)$ to indicate $\ltm(\{ s\})$.

\zorc{
\begin{example}
Let $\lang = \{\wA\wB,\wA\wB\wB\}$ and $\ltm(\wA) = \MA,\ \ltm(\wB)=\MB$, where $\MA,\MB$ are defined as in~\Cref{ex:1}. 
Then, as $\{\wA\wB\wB\}\subseteq \lang$, $\ltm(\wA\wB\wB) = ABB\preceq \ltm(\lang)=AB\oplus ABB$; indeed,
\[
    ABB = \begin{bmatrix} -\infty & -\infty & -\infty \\ 0 & 1 & -\infty \\ -3 & -\infty & -\infty \end{bmatrix}\quad \mbox{and} \quad 
    AB \oplus ABB = \begin{bmatrix} -\infty & -\infty & -\infty \\ 0 & 1 & 3 \\ -3 & -3 & -\infty \end{bmatrix}.
\] 
\end{example}
}

As in precedence graphs, where elements of the $r$-th power of matrices can be interpreted as weights of paths of length $r$, we can interpret elements of products of $r$ matrices as weights of paths of length $r$ in multi--precedence graphs.
For instance, let $ \MA $ and $ \MB $ be $n\times n$ matrices in $\Rmax$ and let us consider the multi--precedence graph $\graph( \MA , \MB )$.
Then, $( \MA   \MB ^2)_{ij} = \max_{k_1,k_2} ( \MA _{ik_1}+ \MB _{k_1k_2}+ \MB _{k_2j})$ is equal to the maximum weight of all paths in the multi--precedence graph $\mathcal{G}( \MA , \MB )$ of the form
\[
\sigma = j \xrightarrow{ \wB } k_2 \xrightarrow{ \wB } k_1 \xrightarrow{ \wA } i
\]
for all $k_1,k_2=1,\ldots,n$.
In the same way, the diagonal element $( \MA  \MB ^2)_{ii}$ represents the maximum weight of all circuits from node $i$ of the form
\[
\sigma = i \xrightarrow{ \wB } k_2 \xrightarrow{ \wB } k_1 \xrightarrow{ \wA } i
\]
for all $k_1,k_2=1,\ldots,n$.
More generally, we can state the following proposition.

\begin{proposition} \label{pr:shared_graphs_properties}
Let $ \MA _1,\ldots, \MA _l$ be $n\times n$ matrices in $\Rmax$ and $\mathcal{G}( \MA _1,\ldots, \MA _l)=(\nodes,\Sigma, \ltm , \arcs)$ the multi--precedence graph associated with them.
Then, for each string $s\in \Sigma^*$, element $i,j$ of matrix $\ltm(s)$,
\begin{equation*}
\begin{split}
\ltm (s)_{ij} & = \left( \ltm (s(1)) \otimes  \ltm (s(2)) \otimes \ldots \otimes  \ltm (s(|s|)) \right)_{ij} = \\ &= \max_{k_1,\ldots,k_{|s|-1}} \left( \ltm (s(1))_{ik_1} +  \ltm (s(2))_{k_1k_2}+\ldots + \ltm (s(|s|))_{k_{|s|-1}j}\right),
\end{split}
\end{equation*}
is equal to the maximum weight of all paths $\sigma$ in $\mathcal{G}( \MA _1,\ldots, \MA _l)$ of the form
\[
\sigma = j \xrightarrow{s(|s|)} k_{|s|-1} \xrightarrow{s(|s|-1)} k_{|s|-2} \xrightarrow{s(|s|-2)} \cdots \xrightarrow{s(2)} k_{1} \xrightarrow{s(1)} i,
\]
for all $k_1,\ldots,k_{|s|-1} \in \{ 1,\ldots,n\}$.
\end{proposition}

\zorc{
As a consequence of the latter proposition and the definition of the morphism $\ltm$ extended to $2^{\Sigma^*}$, given a language $\lang\subseteq \Sigma^*$, 
\[
    \ltm(\lang)_{ij} = \bigoplus_{s\in\lang} \ltm(s)_{ij}
\]
corresponds to the supremum, for all strings $s\in \lang$, of the weights of all paths labeled $s$ in $\graph(A_1,\ldots,A_l)$ from node $j$ to node $i$.
}

The following proposition shows that multi--precedence graphs can be used to study the sign of circuit weights in some precedence graphs.

\begin{proposition} \label{pr:shared_wedge_graphs}
Let $ \MA _1,\ldots, \MA _l$ be $n\times n$ matrices in $\Rmax$ and $\mathcal{G}( \MA _1,\ldots, \MA _l)=(\nodes,\Sigma, \ltm , \arcs_1)$ the multi--precedence graph associated with them.
Then \linebreak$\mathcal{G}( \MA _1,\ldots, \MA _l)$ has a circuit with positive weight from node $i \in \nodes$ iff the precedence graph $\mathcal{G}( \MA _1\oplus\ldots\oplus  \MA _l)=(\nodes,\arcs_2,w)$ has a circuit with positive weight from node $i\in \nodes$.
\end{proposition}
\begin{proof}
The precedence graph $\mathcal{G}( \MA _1\oplus\ldots\oplus  \MA _l)=(\nodes,\arcs_2,w)$ can be built from the multi--precedence graph $\mathcal{G}( \MA _1,\ldots,  \MA _l)=(\nodes,\Sigma, \ltm , \arcs_1)$ as follows:
\begin{itemize}
    \item[-] $\arcs_2 \subseteq \nodes \times \nodes$ is defined such that there is an arc $(j,i)$ from node $j$ to node $i$ iff $\exists \wZ \in \Sigma$ such that $(j,i,\wZ) \in \arcs_1$,
    \item[-] $w:\arcs_2 \rightarrow \R$ is defined for all $(j,i)\in \arcs_2$ as
    \[
    w((j,i)) = \bigoplus_{\wZ\in \Sigma}  \ltm _{ij}(\wZ).
    \]
\end{itemize}
Therefore, if there is a circuit 
\[
\sigma = i \xrightarrow{\wZ_1} k_{r-1} \xrightarrow{\wZ_2} \cdots \xrightarrow{\wZ_{r-1}} k_{1} \xrightarrow{\wZ_{r}} i
\]
in $\mathcal{G}( \MA _1,\ldots,  \MA _l)$ from node $i \in \nodes$ such that $|\sigma|_W>0$, then, by construction, there exists in $\graph( \MA _1\oplus\ldots\oplus  \MA _l)$ a circuit
\[
\rho = i \rightarrow k_{r-1} \rightarrow \cdots \rightarrow k_{1} \rightarrow i
\]
from node $i \in \nodes$, and its weight is such that $|\rho|_W\geq |\sigma|_W >0$.
Conversely, if there is a circuit with positive weight in precedence graph $\mathcal{G}( \MA _1\oplus\ldots\oplus  \MA _l)$, then the same circuit, with the same weight, is present in multi--precedence graph $\mathcal{G}( \MA _1,\ldots,  \MA _l)$.
\end{proof}

Similarly to precedence graphs, we indicate with $\nonegsetm$ the set of multi--precedence graphs that do not contain any circuit with positive weight.

\begin{example}\label{ex:2}
Let us consider matrices $ \MA $, $ \MB $ of Example~\ref{ex:1}.
From  Proposition~\ref{pr:shared_graphs_properties} and
\begin{equation*}
 \MA    \MB  =     
\begin{bmatrix}
-\infty & -\infty & -\infty \\
-\infty & 0 & 3 \\
-\infty & -3 & -\infty
\end{bmatrix},
\end{equation*}
we can conclude, for example, that there is at least a circuit in $\mathcal{G}( \MA , \MB )$ of form
\[
\sigma_1(k) = 2 \xrightarrow{ \wB } k \xrightarrow{ \wA } 2,
\]
$k \in \{1,2,3\}$, and also that the maximum weight for a circuit of this form is $\max_{k} |\sigma_1(k)|_W = 0$.
In particular, the maximum is attained by $\sigma_{1}(1) = 2 \xrightarrow{ \wB } 1 \xrightarrow{ \wA } 2$.
In the same way, since
\begin{equation*}
 \MA    \MB ^2 =     
\begin{bmatrix}
-\infty & -\infty & -\infty \\
0 & 1 & -\infty \\
-3 & -\infty & -\infty
\end{bmatrix},
\end{equation*}
then there exists a circuit in $\mathcal{G}( \MA , \MB )$ of the form
\[
\sigma_2(k_1,k_2) = 2 \xrightarrow{ \wB } k_1 \xrightarrow{ \wB } k_2 \xrightarrow{ \wA } 2,
\]
$k_1,k_2\in\{1,2,3\}$ and the maximum weight for a circuit of this form is \linebreak$\max_{k_1,k_2} |\sigma_2(k_1,k_2)|_W=1$. In particular, the maximum is attained by $\sigma_{2}(3,1)=2 \xrightarrow{ \wB } 3 \xrightarrow{ \wB } 1 \xrightarrow{ \wA } 2$.
Because of the presence of a circuit with positive weight from node $2$ in $\graph(\MA,\MB)$, from  Proposition~\ref{pr:shared_wedge_graphs} we can conclude that there exists a circuit with positive weight from node $2$ in $\graph(\MA\oplus\MB)$, i.e., $(( \MA \oplus  \MB )^*)_{22} = +\infty$.
\end{example}

In the following, we state an interesting property of the Kleene star operator, when applied to a matrix that can be factored as a product of two square matrices.
The proposition will also show how the interpretation of matrix multiplications in terms of weight of paths in multi--precedence graphs can be used to prove algebraic statements in the max-plus algebra.

\begin{proposition} \label{pr:cyclicity}
Let $ \MA, \MB\in \Rmax^{n\times n}$.
Then $\exists i$ such that $((\MA \otimes \MB)^*)_{ii} = +\infty$ if and only if $\exists k$ such that $((\MB \otimes \MA)^*)_{kk} = +\infty$.
\end{proposition}
\begin{proof}
Let us define the multi--precedence graph $\graph( \MA, \MB) = (\nodes, \{ \wA, \wB\}, \ltm ,\arcs)$.
\zorc{
Suppose that $((A\otimes B)^*)_{ii}=+\infty$; then, from~\Cref{pr:aux}, there exists a circuit with positive weight from node $i$ in $\graph(A\otimes B)$, i.e., there exists a number $r\in\mathbb{N}$ such that $((A\otimes B)^r)_{ii}\succ 0$.
\Cref{pr:shared_graphs_properties} implies that there is a circuit with positive weight $\sigma$ in $\graph( \MA, \MB)$
labeled $(\wA\wB)^r$ (i.e., string $\wA\wB$ repeated $r$ times) of the form
\[
\sigma = i \xrightarrow{ \wB} k_1 \xrightarrow{ \wA} k_2 \xrightarrow{ \wB} k_3 \xrightarrow{ \wA} k_4 \xrightarrow{\wB}\dots\xrightarrow{\wB} k_{2r-1} \xrightarrow{\wA}  i
\]
for some $k_1,\dots,k_{2r-1}\in \nodes$.
Let us build another circuit $\tilde{\sigma}$ in $\graph(A,B)$ labeled $(\wB\wA)^r$ as
\[
\tilde{\sigma}  = k_1 \xrightarrow{ \wA} k_2 \xrightarrow{ \wB} k_3 \xrightarrow{ \wA} k_4 \xrightarrow{\wB}\dots\xrightarrow{\wB} k_{2r-1} \xrightarrow{\wA}  i \xrightarrow{\wB} k_1.
\]}
Since the edges in $\sigma$ and $\tilde{\sigma}$ are the same, $|\tilde{\sigma}|_W = |\sigma|_W > 0$; thus \zorc{
\linebreak$(( \MB\otimes  \MA)^r)_{k_1k_1} \succ 0$, which also implies that $(( \MB\otimes  \MA)^*)_{k_1k_1} =+\infty$.
The proof is completed by taking $k=k_1$.} 
\end{proof}

We remark that a consequence of the latter proposition is that, given $\MA,\MB\in\Rmax^{n\times n}$, $\graph(\MA\otimes\MB)\in\nonegset$ if and only if $\graph(\MB\otimes\MA)\in\nonegset$, \zorc{i.e., $tr((\MA \otimes \MB)^*) = tr((\MB\otimes \MA)^*)$}.

\section{Problem formulation}\label{se:main_problem}

Given a parametric precedence graph $\graph(\MA)$, $ \MA  = A(\lambda_1,\ldots,\lambda_p)\in \Rmax^{n\times n}$, the \textit{Non-positive Circuit weight Problem} (NCP) consists in finding the values for $\lambda_1,\ldots,\lambda_p\in \R$ such that $\graph(\MA)\in \nonegset$.
The solution of this problem is already known when $\MA$ depends proportionally ($\MA(\lambda) = \lambda \MA$) or inversely ($\MA(\lambda) = \lambda^{-1} \MA$) on one single parameter, and is recalled as follows.

\begin{proposition}[From Theorem 1.6.18 in~\cite{butkovivc2010max}] \label{pr:minimumcyclemean}
The precedence graph $\mathcal{G}(\lambda \MA)$ belongs to $\nonegset$ iff $\lambda \preceq \left(\mbox{\normalfont mcm} (\MA)\right)^{-1}$.
Similarly, the precedence graph $\mathcal{G}(\lambda^{-1} \MA)$ belongs to $\nonegset$ iff $\lambda \succeq \mbox{\normalfont mcm} (\MA)$.
\end{proposition} 

\zorc{
\begin{remark}\label{re:2}
    It is easy now to prove the following property of the maximum circuit mean: for all $\MA,\MB\in\Rmax^{n\times n}$, 
    \[
    \mbox{mcm}(\MA\otimes \MB) = \mbox{mcm}(\MB\otimes \MA)\ .
    \] 
    Indeed, due to~\Cref{pr:cyclicity}, $\graph(\lambda^{-1}\MA\otimes\MB)\in\nonegset$ iff $\graph(\lambda^{-1}\MB\otimes\MA)\in\nonegset$.
    To complete the proof, we just need to apply~\Cref{pr:minimumcyclemean} to parametric precedence graphs $\graph(\lambda^{-1}\MA\otimes\MB)$ and $\graph(\lambda^{-1}\MB\otimes\MA)$.
\end{remark}
}

In this paper, we are interested in solving the NCP on a class of more general parametric precedence graphs.
Let $ \MP ,\  \MI ,\  \MC $ be three arbitrary $n\times n$ matrices with elements in $\Rmax$.
Let $\graphPIC$ be a parametric precedence graph depending on the only parameter $\lambda\in \R$.
The weights of the arcs in this graph can depend proportionally ($\lambda \MP $), inversely ($\lambda^{-1}  \MI $) and constantly ($ \MC $) with respect to $\lambda$.
For this reason, we refer to the NCP on this class of parametric precedence graphs as \textit{Proportional-Inverse-Constant}-NCP (PIC-NCP).
From Proposition~\ref{pr:shared_wedge_graphs}, the problem can be equivalently stated on the parametric multi--precedence graph $\graphPICm$.
We indicate by $\solPIC$ the set of $\lambda$’s that solve the PIC-NCP on graph $\graphPIC$:
\[
	\solPIC \coloneqq \{\lambda \in\R\ | \ \graphPIC \in\nonegset\} = \{\lambda\in\R \ | \ \graphPICm\in\nonegsetm\}.
\]

\subsection{Solution using linear programming}

In this subsection, we show how to use linear programming in standard algebra to solve the PIC-NCP.

\begin{proposition}
The existence of a solution $\lambda\in \R$ of the PIC-NCP can be checked in standard algebra using linear programming.
\end{proposition}
\begin{proof}
From Proposition~\ref{pr:v_leq_Bv}, $\graphPIC\in \nonegset$ for some $\lambda\in \R$ if and only if there exists a vector $[x^\intercal, \lambda]^\intercal \in \R^{n+1}$ such that
\begin{equation}\label{eq:main_ineq}
x \succeq (\lambda \MP \oplus \lambda^{-1} \MI \oplus \MC) \otimes x.
\end{equation}
The $i$-th inequality can be written in standard algebra as
\[
x_i \geq \max_{j=1,\ldots,n} (\lambda + \MP_{ij} + x_j,-\lambda + \MI_{ij} + x_j,\MC_{ij} + x_j),
\]
which is equivalent to the following linear system
\begin{equation*}
\left\{ 
\begin{array}{ll}
x_i \geq \lambda + \MP_{ij} + x_j & \forall j = 1,\ldots,n\\
x_i \geq -\lambda + \MI_{ij} + x_j & \forall j = 1,\ldots,n\\
x_i \geq \MC_{ij} + x_j & \forall j = 1,\ldots,n
\end{array}
\right.
\end{equation*}
The system above, written for all $i$, forms a system of $m\leq 3n^2$ linear inequalities in $n+1$ variables $x_1,\ldots,x_{n},\lambda$, where $m$ is the number of edges in $\graphPICm$.
Indeed, when an element of the matrices $\MP$, $\MI$ or $\MC$ is $-\infty$, the corresponding inequality is automatically satisfied, as the right hand side becomes $-\infty$.
\end{proof}

An important observation that comes from the latter proposition is that the solution set of Inequality~\eqref{eq:main_ineq} is convex in $[x^\intercal, \lambda]^\intercal$.
Since the projection of a convex set is convex, the consequence is that the set of $\lambda$’s that solves the PIC-NCP is an interval $\solPIC = [\lambda_{\min},\lambda_{\max}]\cap \R$,
where $\lambda_{\min}\in\Rbar$ is the optimal value of the linear program
\begin{align}\label{lp:min}\tag{LP1}
\min \quad & \lambda &\\
\mbox{s.t.}\quad  & x_i \geq \lambda + \MP_{ij} + x_j &\forall i,j=1,\ldots,n\nonumber \\
& x_i \geq -\lambda + \MI_{ij} + x_j   &\forall i,j=1,\ldots,n\nonumber \\
& x_i \geq \MC_{ij} + x_j   &\forall i,j=1,\ldots,n\nonumber 
\end{align}
and, similarly, $\lambda_{\max}\in\Rbar$ is the optimal value of the linear program
\begin{align}\label{lp:max}\tag{LP2}
\max \quad & \lambda &\\
\mbox{s.t.}\quad  & x_i \geq \lambda + \MP_{ij} + x_j &\forall i,j=1,\ldots,n\nonumber \\
& x_i \geq -\lambda + \MI_{ij} + x_j   &\forall i,j=1,\ldots,n\nonumber \\
& x_i \geq \MC_{ij} + x_j   &\forall i,j=1,\ldots,n.\nonumber 
\end{align}
The consequence is that the PIC-NCP can be solved by solving problems~\eqref{lp:min} and~\eqref{lp:max}.
Observe that $\solPIC$ coincides with $[\lambda_{\min},\lambda_{\max}]$, unless $\lambda_{\min}$ or $\lambda_{\max}$ are not finite, in which case $\solPIC$ is an unbounded interval.

\subsection{Solution using dioid theory techniques}

We recall that it is still an open problem to find an algorithm that solves linear programs in strongly polynomial time.
Our aim is to derive a closed-formula expression for $\lambda_{\min}$ and $\lambda_{\max}$ that can be computed in strongly polynomial time, i.e., in a number of operations that does not depend on the size of parameters $\MP_{ij}$, $\MI_{ij}$ and $\MC_{ij}$, for all $i,j$.
The following theorem states that Algorithm~\ref{al:lambdasnn} provides the desired closed formula.
Observe that in line 7 of the algorithm, $\left \lfloor{\cdot } \right \rfloor$ indicates the floor function.

\begin{theorem}\label{th:correctness}
Algorithm~\ref{al:lambdasnn} 
solves the PIC-NCP for graph $\graph(\lambda \MP \oplus \lambda^{-1} \MI  \oplus \MC)$.
\end{theorem}

\newlength{\commentWidth}
\setlength{\commentWidth}{5cm}
\newcommand{\atcp}[1]{\tcp*[r]{\makebox[\commentWidth]{#1\hfill}}}

\begin{algorithm}[t] \DontPrintSemicolon \label{al:lambdasnn} \small
\SetAlgoLined
\KwIn{$\MP,\MI,\MC\in \Rmax^{n\times n}$}
\KwOut{$\Lambda(P, I, C)$}
 \eIf{$\graph(\MC)\notin \nonegset$}{
     return $\emptyset$\\
 }{
    $\MP \leftarrow \MC^* \MP \MC^*$ \atcp{denoted $\ltm ( \xx )$ in Subsection~\ref{su:technical}} 
    $\MI \leftarrow \MC^* \MI \MC^*$ \atcp{denoted $\ltm ( \yy )$ in Subsection~\ref{su:technical}}
    $S \leftarrow E_{\otimes}$       \atcp{initialize $ \ltm ( \xy ( \xx , \yy ,0))$}
    \For{$k = 1 \text{\normalfont{ to }} \left \lfloor{\frac{n}{2}} \right \rfloor$}{\vspace*{2pt}
        $S\leftarrow \MP S^2 \MI \oplus \MI S^2 \MP \oplus E_{\otimes}$ \atcp{compute $ \ltm ( \xy ( \xx , \yy ,k))$}
    }
    \eIf{$\graph(S)\notin \nonegset$}{
     return $\emptyset$
    }{
        $\lambda_{\min} \leftarrow \mbox{mcm}(\MI S^*)$\\
        $\lambda_{\max} \leftarrow (\mbox{mcm}(\MP S^*))^{-1}$\\
        return $[\lambda_{\min},\lambda_{\max}]\cap \R$\\
    }
 }
 \caption{Find all $\lambda$ such that $\graphPIC\in \nonegset$}
\end{algorithm}

The proof will be given in the next section.
In the reminder of this section, we will suppose that the algorithm is correct.
We observe that, given cubic complexity for square matrix max-plus multiplication, the time complexity of the algorithm is $\mathcal{O}(n^4)$, since it is dominated by the ``for'' loop.
Moreover, its space complexity is $\mathcal{O}(n^2)$ as, at each step of the algorithm, $3$ matrices of size $n\times n$ must be stored.

The main idea behind the algorithm is to build a sequence of matrices $M^{(k)}_{\min}$ and $M^{(k)}_{\max}$, such that the sequence of intervals
\[
\pazocal{I}^{(k)}\coloneqq[\mbox{mcm}(M^{(k)}_{\min}),(\mbox{mcm}(M^{(k)}_{\max}))^{-1}]
\]
approximates better and better $[\lambda_{\min},\lambda_{\max}]$ with increasing values of $k\in\nato$.
A first approximation of $[\lambda_{\min},\lambda_{\max}]$ is obtained by setting $M^{(0)}_{\min}\coloneqq \MC^*\MI\MC^*$ and $M^{(0)}_{\max}\coloneqq \MC^*\MP\MC^*$ (lines 4-5)\footnote{These matrices correspond, respectively, to matrices $\pazocal{B}^\sharp$ and $\pazocal{A}$ introduced in~\cite{paek2020analysis}.}. \enlargethispage{-\baselineskip}
Indeed, with this choice we have $[\lambda_{\min},\lambda_{\max}]\subseteq \pazocal{I}^{(0)}\subseteq [\mbox{mcm}(\MI),(\mbox{mcm}(\MP))^{-1}]$. 
In order to get tighter approximations, it is necessary to introduce another sequence of matrices: let $S^{(k)}$ denote the matrix obtained after $k\in\nato$ iterations of the "for" loop of lines 7-9, with $S^{(0)}=E_\otimes$.
By choosing, for all $k\in\nato$, $M^{(k)}_{\min}\coloneqq M^{(0)}_{\min}(S^{(k)})^*$ and $M^{(k)}_{\max}\coloneqq M^{(0)}_{\max}(S^{(k)})^*$, we have $[\lambda_{\min},\lambda_{\max}]\subseteq\pazocal{I}^{(k+1)}\subseteq \pazocal{I}^{(k)}$.
Moreover, as it will be proved in the next section, the sequence of intervals $\pazocal{I}^{(k)}$ converges to $[\lambda_{\min},\lambda_{\max}]$ after at most $k=\left \lfloor{\frac{n}{2}} \right \rfloor$ iterations.

\begin{remark}
In case $\MP$, $\MI$ or $\MC$ is $\pazocal{E}$, the problem reduces to a subclass of the PIC-NCP.
Clearly, if Algorithm~\ref{al:lambdasnn} solves the PIC-NCP, then it solves these problems, too.
Moreover, in these cases, the algorithm can be simplified.
In particular, it can be proven that\footnote{
\zorc{
    To do so, recall~\Cref{re:2} and~\eqref{eq:prop0}.
}}:
\begin{itemize}
\item[-] if $\MP = \pazocal{E}$, the complexity of the algorithm reduces to $\pazocal{O}(n^3)$, since
    \[
\begin{array}{ll}
\lambda_{\min} = \mbox{mcm}(\MI\MC^*), &  \quad \lambda_{\max} = +\infty,
\end{array}
\]
\item[-] if $\MI = \pazocal{E}$, the complexity of the algorithm reduces to $\pazocal{O}(n^3)$, since
\[
\begin{array}{ll}
\lambda_{\min} = -\infty, & \quad  \lambda_{\max} = (\mbox{mcm}(\MP\MC^*))^{-1},
\end{array}
\]
\item[-] if $\MC = \pazocal{E}$, lines 1-5 of the algorithm are no longer needed, but its complexity remains the same.
\end{itemize}

Observe that, as Levner-Kats algorithm can be used to solve the generic PIC-NCP at the cost of increasing its worst-time complexity, the two algorithms from~\cite{karp1981parametric} and the one from~\cite{young1991faster} can be used to solve the PIC-NCP when $\MP = \pazocal{E}$ or $\MI = \pazocal{E}$.
However, nodes and arcs must be added to the original graph, as in Figure~\ref{fig:LevnerKatsa}.
This increases the worst-case complexity of the algorithms, respectively, to $\pazocal{O}(n^6)$, $\pazocal{O}(n^4\log n)$ and $\pazocal{O}(n^4\log n)$.
Therefore, our algorithm solves also these subclasses of the PIC-NCP faster than traditional ones, in the worst-case.
\end{remark}

\begin{example}\label{ex:3}
Let
\[
\MP = \begin{bmatrix}
     -\infty & -\infty & -\infty\\
     -\infty & -\infty & -\infty\\
     -\infty & -\infty & -4
    \end{bmatrix},\ 
\MI = \begin{bmatrix}
     -\infty & 0 & -\infty\\
     -\infty & -\infty & 0.5\\
     -\infty & -\infty & 0
    \end{bmatrix},\ 
\MC = \begin{bmatrix}
    -\infty & -3 & -\infty\\
    2 & -\infty & -\infty\\
    6 & 0.5 & -\infty
    \end{bmatrix}.
\] 
The aim is to find all the values of $\lambda$ such that the parametric precedence graph $\graphPIC$ does not contain circuits with positive weight.
From~\Cref{pr:shared_wedge_graphs}, we can do so by studying the parametric multi--precedence graph $\graph(\lambda \MP,\lambda^{-1} \MI,\MC) = (\nodes,\{\wP,\wI,\wC\}, \ltm , \arcs) $, which is depicted in Figure~\ref{fig:ex3}, where $ \ltm (\wP) = \lambda \MP$, $ \ltm (\wI) = \lambda^{-1} \MI$ and $ \ltm (\wC) = \MC$.
\begin{figure}
\centering
\resizebox{.6\linewidth}{!}{
\begin{tikzpicture}[node distance=1cm and 2cm,>=stealth',bend angle=45,thick]

\node [place,tokens=0,label=center:{$1$},] (P1) {};
\node [place,tokens=0,label=center:{$2$},right= of P1] (P2) {};
\node [place,tokens=0,label=center:{$3$},right= of P2] (P3) {};

\draw (P2) edge[->] node[auto,swap,label=below:{$-\lambda,\wI$}] {} (P1);
\draw (P3) edge[loop below,->] node {$-\lambda,\wI$} (P3);
\draw (P3) edge[loop right,->] node {$-4+\lambda,\wP$} (P3);
\draw (P2) edge[bend left=60,->] node[auto,swap,label=below:{$-3,\wC$}] {} (P1);
\draw (P1) edge[bend left=60,->] node[auto,swap,label=above:{$2,\wC$}] {} (P2);
\draw (P3) edge[->] node[auto,swap,label=below:{$0.5-\lambda,\wI$}] {} (P2);
\draw (P1) edge[bend left=90,->] node[auto,swap,label=above:{$6,\wC$}] {} (P3);
\draw (P2) edge[bend left=60,->] node[auto,swap,label=above:{$0.5,\wC$}] {} (P3);

\end{tikzpicture}
}
\caption{Parametric multi--precedence graph $\graph(\lambda \MP,\lambda^{-1} \MI,\MC)$ from Example~\ref{ex:3}.}
\label{fig:ex3}
\end{figure}
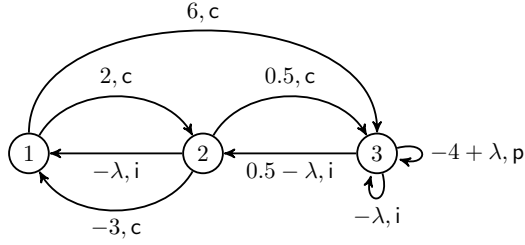
The interval $\solPIC = [\lambda_{\min},\lambda_{\max}] \cap \R$
can be obtained by using Algorithm~\ref{al:lambdasnn}.
We get
\[
\lambda_{\min} = \mbox{mcm}(\tilde{\MI}S^*) \quad \lambda_{\max} = (\mbox{mcm}(\tilde{\MP}S^*))^{-1},
\]
where
\[
\tilde{\MI} = \MC^*\MI\MC^* =\begin{bmatrix}
3.5 & 0.5 & -2.5\\
6.5 & 3.5 & 0.5\\
9.5 & 6.5 & 3.5
\end{bmatrix},\
\tilde{\MP} = \MC^*\MP\MC^* =\begin{bmatrix}
-\infty & -\infty & -\infty\\
-\infty & -\infty & -\infty\\
2 & -1 & -4
\end{bmatrix},
\]
\[
S^* = (\tilde{\MP}\tilde{\MI} \oplus \tilde{\MI}\tilde{\MP} \oplus E_{\otimes})^* = \begin{bmatrix}
0 & -3.5 & -6.5\\
2.5 & 0 & -3.5\\
5.5 & 2.5 & 0
\end{bmatrix};
\]
hence $\lambda_{\min} = 3.5$, $\lambda_{\max} = 4$.
\end{example}

\section{Correctness of the algorithm}\label{se:validity}

In this section, we prove Theorem~\ref{th:correctness}.
The proof is based on the connection between formal languages and multi--precedence graphs.
The propositions are divided in propositions on formal languages (Subsection~\ref{su:formallanguages}) and propositions on graphs (Subsection~\ref{su:technical}).

\zorc{
In the following, we give a sketch of the proof, hoping it will help the reader delve into the technical details contained in the next subsections.
To simplify the discussion, here we consider the case in which $\MC=\pazocal{E}$.
In Subsection~\ref{su:formallanguages} we will prove that, given a binary alphabet $\Sigma=\{\wA,\wB\}$, a number $n\in\nato$, and a language $\xy$ containing only balanced strings and such that $\xy^*$ contains all balanced strings from $\Sigma^*$ of length $\leq n$ (we will give an example of such language in~\Cref{def:S}),
\[
\sum_{k=0}^{n} (\wA+\wB)^k \subseteq \xy^* + (\xy^*\wA \xy^*)^* + (\xy^*\wB \xy^*)^* \subseteq (\wA+\wB)^*,
\] 
i.e., any word from $\Sigma^*$ of length $\leq n$ either belongs to $\xy^*$ (if it is balanced), or to $(\xy^*\wA \xy^*)^*$ (if it contains more $\wA$'s than $\wB$'s) or to $(\xy^*\wB \xy^*)^*$ (if it contains more $\wB $'s than $\wA$'s).
In Subsection~\ref{su:technical}, by applying morphism $\ltm$ to the latter expression, with $\ltm(\wA)=\lambda\MP$, $\ltm(\wB)=\lambda^{-1}\MI$, $\MP,\MI\in\Rmax^{n\times n}$, and $\ltm(\xy)=S$ (matrix $S$ obtained in this way will coincide with the homonym matrix computed in lines 6-9 of Algorithm~\ref{al:lambdasnn}) and using~\Cref{pr:cyclicity} together with the fact that $A\preceq B \Rightarrow tr(A)\preceq tr(B)$, we will show that
\begin{align*}
    tr \left( \bigoplus_{k=0}^{n} (\lambda\MP\oplus \lambda^{-1}\MI)^k \right) & \preceq tr \left( S^*\oplus (S^*\lambda\MP S^*)^*\oplus (S^*\lambda^{-1}\MI S^*)^* \right) = \\
& = \bigoplus_{i=1}^{n} \underbrace{ \left( (S^*)_{ii}\oplus ((\lambda\MP S^*)^*)_{ii}\oplus ((\lambda^{-1}\MI S^*)^*)_{ii} \right) }_{\coloneqq M_{ii}} \preceq \\
& \preceq tr \left( (\lambda\MP\oplus\lambda^{-1}\MI)^* \right) .
\end{align*}
Note that, on the one hand, if $M_{ii}=0$ for all $i$, then there are no (elementary) circuits of positive weight in $\graph(\lambda\MP\oplus\lambda^{-1}\MI)$, as $tr \left( \bigoplus_{k=0}^{n} (\lambda\MP\oplus \lambda^{-1}\MI)^k \right) \preceq tr(M)$;
on the other hand, if there exists $i$ such that $M_{ii}=+\infty$, then there are circuits with positive weight in $\graph(\lambda\MP\oplus\lambda^{-1}\MI)$, as $tr(M) \preceq tr\left( (\lambda\MP\oplus\lambda^{-1}\MI)^*\right)$.
Therefore, $\graph(\lambda\MP\oplus\lambda^{-1}\MI)\in\nonegset$ if and only if $tr(M)=0$.
Moreover, since $\xy$ contains only balanced strings, matrix $S = \ltm(\xy)$ does not depend on $\lambda$, as for each string $s\in\xy$, $\lambda$ and $\lambda^{-1}$ cancel out in $\ltm(s)$ (as $\lambda\otimes \lambda^{-1}=\lambda^{-1}\otimes \lambda=0$).
Hence, from~\Cref{pr:minimumcyclemean}, $tr(M)=0$ if and only if $\graph(S)\in\nonegset$ and $\mbox{mcm}(\MI S^*)\preceq\lambda\preceq(\mbox{mcm}(\MP S^*))^{-1}$.
} 

\subsection{Propositions on formal languages} \label{su:formallanguages}

In the fi\zorc{r}st part of this subsection, we focus on balanced binary strings.

\begin{lemma} \label{le:balStrings}
Every $x*x$--balanced binary string of positive length can be built by concatenating two $x*y$--balanced binary strings of positive length.
\end{lemma}
\begin{proof}

Let $\Sigma = \{\wA,\wB\}$ be a binary alphabet.
Let us define a function $g:\Sigma \rightarrow \{-1,+1\}$ as
\[
g(\wZ) = \begin{cases}
       	+1 & \mbox{ if } \wZ = \wA \\
       	-1 & \mbox{ if } \wZ = \wB,
       \end{cases}
\]
and function $h_s:\{1,\ldots,|s|\} \rightarrow \mathbb{Z}$, associated with a string $s\in \Sigma^*\setminus\{\emptystr\}$, as
\[
h_s(i) = \sum_{j = 1}^i g(s(j)).
\]
Intuitively, $h_s$ counts up or down by $1$ if $g(s(i))$ is $+1$ or $-1$.
In Figure~\ref{fig:wordGraph}, $h_s$ is plotted for a given string $s$.
It is clear that a binary string $s \in \Sigma^*\setminus\{\emptystr\}$ is balanced iff $h_s(|s|) = 0$.
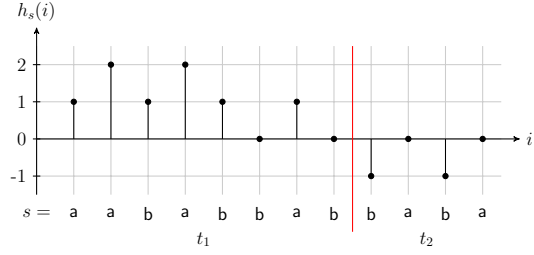
\begin{figure}
\centering
\resizebox{.6\linewidth}{!}{
\begin{tikzpicture}[node distance=.5cm and 1cm,>=stealth',bend angle=45,thick]
\newcommand*\circled[1]{\tikz[baseline=(char.base)]{
            \node[shape=circle,draw,inner sep=2pt,red] (char) {\textcolor{black}{#1}};}}
\Large{
\draw[step=1,gray!40](0,-1.5) grid (12.5,2.5); 
\draw[->] (0,0) -- (13,0) node[right] {$i$};
\draw[->] (0,-1.5) -- (0,3) node[above] {$h_s(i)$};

\node at (0,-2) {$s=$};
\node at (1,-2) {$\wA$};
\node at (2,-2) {$\wA$};
\node at (3,-2) {$\wB$};
\node at (4,-2) {$\wA$};
\node at (5,-2) {$\wB$};
\node at (6,-2) {$\wB$};
\node at (7,-2) {$\wA$};
\node at (8,-2) {$\wB$};
\node at (9,-2) {$\wB$};
\node at (10,-2) {$\wA$};
\node at (11,-2) {$\wB$};
\node at (12,-2) {$\wA$};

\draw[red] (8.5,-2.5) -- (8.5,2.5);
\node at (4.5,-2.7) {$t_1$};
\node at (10.5,-2.7) {$t_2$};

\filldraw (1,1) circle (2pt);
\filldraw (2,2) circle (2pt);
\filldraw (3,1) circle (2pt);
\filldraw (4,2) circle (2pt);
\filldraw (5,1) circle (2pt);
\filldraw (6,0) circle (2pt);
\filldraw (7,1) circle (2pt);
\filldraw (8,0) circle (2pt);
\filldraw (9,-1) circle (2pt);
\filldraw (10,0) circle (2pt);
\filldraw (11,-1) circle (2pt);
\filldraw (12,0) circle (2pt);

\draw (1,0) -- (1,1);
\draw (2,0) -- (2,2);
\draw (3,0) -- (3,1);
\draw (4,0) -- (4,2);
\draw (5,0) -- (5,1);
\draw (7,0) -- (7,1);
\draw (9,0) -- (9,-1);
\draw (11,0) -- (11,-1);

\foreach \y in {-1,...,2}
     		\draw (0.1,\y) -- (-0.1,\y)
			node[anchor=east] {\y};

}
\end{tikzpicture}
}
\caption{Plot of $h_s$ associated with an $x*x$--balanced string $s=\wA\wA\wB\wA \wB\wB\wA \wB\wB\wA \wB\wA$. The red line divides the string $s$ in two $x*y$--balanced substrings $t_1=\wA\wA\wB\wA \wB\wB\wA \wB$ and $t_2=\wB\wA \wB\wA$.}
\label{fig:wordGraph} 
\end{figure}
Let $s$ be an $x*x$--balanced binary string such that $s(1) = s(|s|) = \wA$, the case in which $s(1) = s(|s|) = \wB$ is analogous.
Then
\begin{equation} \label{eq:h_s_1_len_s}
h_s(1) = +1, \quad h_s(|s|-1) = -1,
\end{equation}
and this implies that $\exists i \in \{2,\ldots,|s|-1\}$ such that
\begin{equation} \label{eq:h_s_0}
	h_s(i-1) = +1 \quad h_s(i) = 0 \quad h_s(i+1) = -1.
\end{equation}
Indeed, since $\forall i\in \{1,\ldots,|s|-1\}$, $|h_s(i+1)-h_s(i)|= 1$, there must be an $i$ that satisfies~\eqref{eq:h_s_0} in order to change the sign of $h_s$ from positive (in $h_s(1)$) to negative (in $h_s(|s|-1)$).
Now, let us define $t_1 = s(1) \ldots s(i)$ and $t_2 = s(i+1) \ldots s(|s|)$. Both $t_1$ and $t_2$ are balanced because $h_{t_1}(|t_1|)=0$ and $h_{t_2}(|t_2|)=0$; moreover $t_1(1)=\wA \neq \wB = t_1(|t_1|)$ and $t_2(1)= \wB \neq \wA = t_2(|t_2|)$.
\end{proof}

\begin{definition}\label{def:S}
Given an alphabet $\Sigma$ and two languages $\lang_1,\lang_2\subseteq \Sigma^*$, we define the sequence of languages $ \xy :2^{\Sigma^*}\times 2^{\Sigma^*} \times \nato \rightarrow 2^{\Sigma^*}$ recursively as
\[
 \xy (\lang_1,\lang_2,0) = \{\emptystr\},
\]
\[
\xy (\lang_1,\lang_2,k) = \lang_1 \xy (\lang_1,\lang_2,k-1)^2\lang_2 + \lang_2 \xy (\lang_1,\lang_2,k-1)^2\lang_1 + \emptystr \quad \forall k\in \nat.
\]
\end{definition}

\zorc{
For instance, for $\lang_1=\{\wA\},\ \lang_2=\{\wB\}$, the first three terms of the sequence are
\[
\xy(\lang_1,\lang_2,0)=\{\emptystr\},\quad \xy(\lang_1,\lang_2,1)=\xy(\lang_1,\lang_2,0)\cup\{\wA\wB,\wB\wA\},  
\] 
\begin{align*}
    \xy(\lang_1,\lang_2,2)= \xy(\lang_1,\lang_2,1)\cup \{& 
\wA\wA\wB\wB,\wA\wB\wA\wB,\wA\wA\wB\wA\wB\wB,\wA\wA\wB\wB\wA\wB,\wA\wB\wA\wA\wB\wB,\wA\wB\wA\wB\wA\wB,\\
&\wB\wA\wB\wA,\wB\wB\wA\wA,\wB\wA\wB\wA\wB\wA,\wB\wA\wB\wB\wA\wA,\wB\wB\wA\wA\wB\wA,\wB\wB\wA\wB\wA\wA\}.
\end{align*}
} 
\begin{theorem} \label{th:s_n}
Given $\Sigma = \{\wA,\wB\}$, $k\in \nato$, language $ \xy (\{\wA\},\{\wB\},k)$ contains all $x*y$--balanced binary strings in $\Sigma^*$ of length less than or equal to $2k$.
\end{theorem}
\begin{proof}
We will use the notation $ \xy _k$ in place of $ \xy (\{\wA\},\{\wB\},k)$ for all $k \in \nato$.
The theorem will be proven by induction.
For $k=0$ the proof is trivial as the empty string, by definition, is $x*y$--balanced.
Suppose that $ \xy _k$ contains all $x*y$--balanced binary strings of length $\leq 2k$, we want to prove that $ \xy _{k+1}$ contains all $x*y$--balanced binary strings of length $\leq 2(k+1)$.

We can write 
\begin{align*}\label{eq:s_n1}
 \xy _{k+1} &= \wA  \xy _k^2 \wB + \wB  \xy _k^2 \wA +\emptystr = \nonumber \\
&= \wA(\wA  \xy _{k-1}^2 \wB + \wB  \xy _{k-1}^2 \wA + \emptystr)^2\wB + \wB(\wA  \xy _{k-1}^2 \wB + \wB  \xy _{k-1}^2 \wA + \emptystr)^2 \wA + \emptystr = \nonumber \\
&= \wA((\wA  \xy _{k-1}^2 \wB + \wB  \xy _{k-1}^2 \wA)^2 + \wA  \xy _{k-1}^2 \wB + \wB  \xy _{k-1}^2 \wA + \emptystr) \wB + \\
&+ \wB((\wA \xy _{k-1}^2\wB+ \wB  \xy _{k-1}^2 \wA)^2 + \wA  \xy _{k-1}^2 \wB + \wB  \xy _{k-1}^2 \wA + \emptystr) \wA + \emptystr = \nonumber \\
&= \wA(\lang^2 + \lang + \emptystr)\wB + \wB(\lang^2 + \lang + \emptystr)\wA + \emptystr \nonumber , 
\end{align*}
where we used the substitution $\lang \coloneqq \wA  \xy _{k-1}^2 \wB + \wB  \xy _{k-1}^2 \wA$.
Since every $x*y$--balanced string of length $\leq 2(k+1)$ either starts with $\wA$ and ends with $\wB$ or starts with $\wB$ and ends with $\wA$, if we prove that $ \xy _k^2 = \lang^2 + \lang + \emptystr$ contains all balanced strings of length $\leq 2k$ (without any condition on the first and last letters) then the proof of the theorem is completed.

Since $\lang =  \xy _k \setminus \{\emptystr \}$, then $\lang$ contains all $x*y$--balanced strings of positive length $\leq 2k$.
From Lemma~\ref{le:balStrings}, every $x*x$--balanced string of positive length can be built by concatenating two $x*y$--balanced strings of positive length; therefore $\lang^2 + \emptystr= \lang\lang + \emptystr$ contains all $x*x$--balanced strings of length $\leq 2k$.
Then, $\lang^2 + \lang + \emptystr$ contains all balanced strings of length $\leq 2k$.
\end{proof}

In the proof, we also showed that $ \xy (\{\wA\},\{\wB\},k)^2$ contains all balanced strings of $\Sigma$ of length $\leq 2k$.
Since $ \xy (\{\wA\},\{\wB\},k)^2 \subseteq  \xy (\{\wA\},\{\wB\},k)^*$, this is also valid for $ \xy (\{\wA\},\{\wB\},k)^*$.
Moreover, since the length of every balanced binary string is an even number, then, given $n$ such that $\left\lfloor \frac{n}{2} \right\rfloor=k$, $ \xy (\{\wA\},\{\wB\},k)$ contains all the $x*y$--balanced strings of length $\leq n$ and the following statement holds.

\begin{corollary}\label{co:all_bal}
Given $\Sigma = \{\wA,\wB\}$, $n\in \nato$, language $ \xy (\{\wA\},\{\wB\},{\left \lfloor{\frac{n}{2}} \right \rfloor})^*$ contains all balanced binary strings in $\Sigma^*$ of length less than or equal to $n$.
\end{corollary}


The second part of this subsection is concerned with "unbalanced" binary strings.

\begin{lemma} \label{pr:strings_division}
Let $\Sigma = \{ \wA ,\wB\}$.
Every binary string $s \in \Sigma^*$ such that $|s|_{ \wA } > |s|_{\wB}$ can be written as
\begin{equation*} \label{eq:strings_decomposition}
	s = t_1  \wA t_2  \wA \cdots  \wA t_r,
\end{equation*}
where $r \in \nat$ and $t_1,\ldots,t_r \in \Sigma^*$ are balanced binary strings.
\end{lemma}

\begin{proof}
Given a binary string $s$ such that $|s|_{ \wA } > |s|_{\wB}$, we define the set $H_s$ as
\[
	H_s\coloneqq \{h\in \{1,2,\ldots,|s|\}\ | \ |s(1)\cdots s(h)|_\wA - |s(1)\cdots s(h)|_\wB = 1 \}.
\]
Note that, since $|s|_{ \wA } - |s|_{\wB}\geq 1$, the set $H_s$ must be non-empty.
Let $h_{\min} \coloneqq \min_h H_s$; then $s(h_{\min}) = \wA$, because if $s(h_{\min}) = \wB$ then 
\[
|s(1)\cdots s(h_{\min}-1)|_\wA-|s(1)\cdots s(h_{\min}-1)|_\wB = 2,
\]
which implies that $\exists h < h_{\min}$ such that $h\in H_s$.
Therefore, it is always possible to factorize the string $s$ as $s = s_1 \wA s_2$, where $s_1 = s(1) \cdots s(h_{\min}-1)$ if $h_{\min} \neq 1$, $s_1 = \emptystr$ otherwise, and $s_2 = s(h_{\min}+1) \cdots s(|s|)$ if $h_{\min}\neq |s|$, $s_2 = \emptystr$ otherwise.
Moreover, since $h_{\min}\in H_s$ and $s(h_{\min})=\wA$, 
\begin{equation*}\label{eq:string_factorization_aux}
	|s_1|_\wA - |s_1|_\wB = 0 \quad \mbox{and} \quad |s_2|_\wA - |s_2|_\wB = |s|_\wA - |s|_\wB - 1.
\end{equation*}
%
Let $t_1 = s_1$ and $m=|s|_\wA-|s|_\wB$.
We will prove the proposition by induction on $m$.

For the case $m = 1$, let $t_2 = s_2$; we have $s = t_1 \wA t_2$ where $t_1$ and $t_2$ are balanced.
Now suppose that the proposition holds for all words such that $|s|_{\wA} - |s|_{\wB} = m$; let us prove that it holds for all words such that $|s|_{\wA} - |s|_{\wB} = m+1$.
As $s = t_1 \wA s_2$ and $t_1$ is balanced, we only need to prove that $s_2$ can be written as $s_2 = t_2 \wA t_3 \wA \cdots \wA t_r$, where $t_2,\ldots,t_r$ are balanced strings; note that this is a direct consequence of the induction hypothesis, since $|s_2|_{\wA} - |s_2|_{\wB} = m$.
\end{proof}

\begin{theorem}\label{pr:main_strings}
Given $\Sigma = \{\wA,\wB\}$, $n\in \nato$, $\xy = \xy (\{\wA\},\{\wB\},{\left \lfloor{\frac{n}{2}} \right \rfloor})$, language $(\xy^*\wA\xy^*)^*$ contains all binary strings $s\in\Sigma^*$ of length less than or equal to $n$ such that $|s|_\wA > |s|_\wB$.
\end{theorem}
\begin{proof}
From~\Cref{pr:strings_division}, any word $s$ with length $n$ or less and such that $|s|_\wA > |s|_\wB$ can be factorized as $s = t_1 \wA t_2 \wA \cdots \wA t_r$, where $t_1,\ldots,t_r$ are balanced strings.
Of course, the length specification imposes that $t_1,\ldots,t_r$ have all length less than or equal to $n$.
Therefore, from~\Cref{co:all_bal}, they must belong to $\xy^*$.
Since 
\begin{align*}
(\xy^*\wA\xy^*)^* &= \emptystr + \xy^*\wA\xy^* + \xy^*\wA\xy^*\xy^*\wA\xy^*+\ldots = \mbox{[from~\eqref{eq:prop0}]} = \\
&=\emptystr + \xy^*\wA\xy^* + \xy^*\wA\xy^*\wA\xy^*+\ldots
\end{align*}
contains all strings formed by concatenating a finite number of time any word from $\xy^*$ with letter $\wA$, and terminating with a word from $\xy^*$, the consequence is that $s=t_1 \wA t_2 \wA \cdots \wA t_r\in(\xy^*\wA\xy^*)^*$.
\end{proof}


%
%

\subsection{Propositions on graphs}\label{su:technical}

In this subsection, we conclude the proof of the correctness of Algorithm~\ref{al:lambdasnn}.
For the reminder of the section, $\MP$, $\MI$ and $\MC$ indicate three arbitrary $n\times n$ matrices in $\Rmax$.

\begin{lemma}\label{le:extension2}
$\graph(\MP \oplus \MI \oplus \MC)\in\nonegset$ if and only if $\graph(\MC)\in\nonegset$ and $\graph(\MC^*\MP\MC^*\oplus \MC^*\MI\MC^*)\in\nonegset$.
\end{lemma}
\begin{proof}
In the following, we will show that $(\MP\oplus \MI\oplus \MC)^* = \MC^* \oplus (\MC^*\MP\MC^* \oplus \MC^*\MI\MC^*)^*$.
From~\Cref{pr:aux}, this will be sufficient for proving the lemma.
Observe that 
\begin{align}\label{eq:sigmastar}
(\MP \oplus\MI \oplus\MC )^* &= \mbox{[from~\eqref{eq:prop1}]} = ( \MC ^*(\MP \oplus\MI ))^*\MC^* =\nonumber\\
&= \mbox{[from~\eqref{eq:prop2}]} = \MC^*((\MP \oplus\MI )\MC^*)^* = \nonumber \\ 
& =\mbox{[from~\eqref{eq:prop3}]} = \MC^*\oplus\MC^*(\MP \oplus\MI)(\MP \oplus\MI \oplus\MC)^*; 
\end{align}
moreover,
\begin{align*}
(\MC ^*\MP \MC ^*\oplus\MC ^*\MI \MC ^*)^* &= (\MC ^*(\MP \oplus\MI )\MC ^*)^* = \mbox{[from~\eqref{eq:prop3}]} =\\ &=E_\otimes \oplus \MC ^*(\MP \oplus\MI )(\MC ^*(\MP \oplus\MI )\oplus\MC )^* = \\
&= E_\otimes \oplus\MC ^*(\MP \oplus\MI )(\MP \oplus\MI \oplus\MC  \oplus \MC ^+(\MP \oplus\MI ))^* = \\
&= \text{[since }(\MP\oplus\MI\oplus\MC)^*=(\MP \oplus\MI \oplus\MC  \oplus \MC ^+(\MP \oplus\MI ))^*\text{]}=\\
&=E_\otimes \oplus \MC ^*(\MP \oplus\MI )(\MP \oplus\MI \oplus\MC )^*.
\end{align*}
Therefore, Equation~\eqref{eq:sigmastar} can be rewritten as
\begin{align*}
(\MP \oplus\MI \oplus\MC )^* &= \MC ^*\oplus\MC ^*(\MP \oplus\MI )(\MP \oplus\MI \oplus\MC )^* =\\
&= \MC ^* \oplus E_\otimes \oplus \MC ^*(\MP \oplus\MI )(\MP \oplus\MI \oplus\MC )^* = \\
&= \MC ^* \oplus (\MC ^*\MP \MC ^*\oplus\MC ^*\MI \MC ^*)^*.
\end{align*}
\end{proof}

\begin{lemma}\label{le:extension3}
Let $\mathcal{G}(\MP,\MI) = (\nodes,\{\wP ,\wI \}, \ltm , \arcs)$, $ \xy  =  \xy ( \{\wP\} , \{\wI\} ,{\left \lfloor{\frac{n}{2}} \right \rfloor})$.
Then, $\graph(\MP\oplus\MI)\in \nonegset$ if and only if $\graph(\ltm(\xy))\in\nonegset$, $\graph(\ltm(\wP\xy^*))\in\nonegset$ and $\graph(\ltm(\wI\xy^*))\in\nonegset$.
\end{lemma}
\begin{proof}
"$\Rightarrow$": since $\xy^*$, $(\wP\xy^*)^*$ and $(\wI\xy^*)^*$ are subsets of $\{\wP,\wI\}^* = (\wP+\wI)^*$, then $\ltm(\xy)^*\preceq (\MP\oplus\MI)^*$, $\ltm(\wP \xy^*)^*\preceq (\MP\oplus\MI)^*$ and $\ltm(\wI \xy^*)^*\preceq (\MP\oplus\MI)^*$.
Therefore, if there exists $i$ such that $(\ltm(\xy)^*)_{ii}=+\infty$, $(\ltm(\wP \xy^*)^*)_{ii}=+\infty$, or $(\ltm(\wI \xy^*)^*)_{ii}=+\infty$, then $((\MP\oplus\MI)^*)_{ii} = +\infty$; from~\Cref{pr:aux} this is equivalent to: $\graph(\ltm(\xy))\notin\nonegset$, $\graph(\ltm(\wP\xy^*))\notin\nonegset$ or $\graph(\ltm(\wI\xy^*))\notin\nonegset$ implies $\graph(\MP\oplus\MI)\notin \nonegset$. 

"$\Leftarrow$": from~\Cref{co:all_bal}, $\xy^*$ contains all the balanced strings of length $\leq n$. From~\Cref{pr:main_strings}, $(\xy^*\wP\xy^*)^*$ contains all the strings of length $\leq n$ with more $\wP$’s than $\wI$’s, and $(\xy^*\wI\xy^*)^*$ contains all the strings of length $\leq n$ with more $\wI$’s than $\wP$’s.
Therefore, for every elementary circuit $\sigma$ from any node $i$ of the multi--precedence graph $\graph(\MP,\MI)$, if its label $s$ is such that $|s|_\wP=|s|_\wI$, then $|\sigma|_W\preceq \ltm(s)_{ii}\preceq (\ltm(\xy)^*)_{ii}=0$, else if $|s|_\wP>|s|_\wI$ then $|\sigma|_W\preceq \ltm(s)_{ii}\preceq (\ltm(\xy^*\wP\xy^*)^*)_{ii}=0$, otherwise $|\sigma|_W\preceq \ltm(s)_{ii}\preceq (\ltm(\xy^*\wI\xy^*)^*)_{ii}=0$.
Finally, from~\Cref{pr:cyclicity}, note that $\graph(\ltm(\xy^*\wP\xy^*))\in\nonegset$ if and only if $\graph(\ltm(\wP\xy^*\xy^*)) = \graph(\ltm(\wP\xy^*))\in\nonegset$, and, similarly, $\graph(\ltm(\xy^*\wI\xy^*))\in\nonegset$ if and only if $\graph(\ltm(\wI\xy^*))\in\nonegset$.
\end{proof}

\begin{theorem} \label{th:extension4}
Let $\mathcal{G}(\MP,\MI,\MC) = (\nodes,\{\wP ,\wI ,\wC \}, \ltm , \arcs)$, $ \xx =\wC ^*\wP \wC ^*$, $ \yy =\wC ^*\wI \wC ^*$, $ \xy  =  \xy ( \xx , \yy ,{\left \lfloor{\frac{n}{2}} \right \rfloor})$, and $\lambda\in\R$.
Then, $\graphPICm\in \nonegsetm$ if and only if $\graph(\MC)\in\nonegset$, $\graph(\ltm(\xy))\in\nonegset$ and $\lambda \in [\lambda_{\min},\lambda_{\max}]\cap \R$, where
\[
\lambda_{\min} = \mbox{\normalfont mcm}( \ltm (\yy \xy^*)) \quad \mbox{ and } \quad \lambda_{\max} = (\mbox{\normalfont mcm}( \ltm (\xx \xy^* )))^{-1}.
\]
\end{theorem}
\begin{proof}
Let $\graphPICm =(\nodes,\{\wP _\lambda,\wI _\lambda,\wC \}, \ltm _\lambda, \arcs)$, where $ \ltm _\lambda(\wP _\lambda) = \lambda \MP$, $ \ltm _\lambda(\wI _\lambda) = \lambda^{-1} \MI $, $ \ltm _\lambda(\wC ) = \MC $.
Moreover, let $ \xx _\lambda=\wC ^*\wP _\lambda \wC ^*$, $ \yy _\lambda=\wC ^*\wI _\lambda \wC ^*$ and \linebreak$ \xy _\lambda =  \xy ( \xx _\lambda, \yy _\lambda,{\left \lfloor{\frac{n}{2}} \right \rfloor})$.
From Lemma~\ref{le:extension2}, $\graphPICm \in \nonegsetm$ if and only if $\graph(\MC )\in \nonegset$ and $\graph(\MC^*\lambda\MP\MC^*\oplus\MC^*\lambda^{-1}\MI\MC^* )\in \nonegset$.
Moreover, from~\Cref{le:extension3}, $\graph(\MC^*\lambda\MP\MC^*\oplus\MC^*\lambda^{-1}\MI\MC^* )=\graph(\ltm_\lambda(\xx_\lambda)\oplus\ltm_\lambda(\yy_\lambda))\in \nonegset$ if and only if $\graph(\ltm_\lambda(\xy_\lambda))\in\nonegset$, $\graph(\ltm_\lambda(\xx_\lambda\xy_\lambda^*))\in\nonegset$ and $\graph(\ltm_\lambda(\yy_\lambda\xy_\lambda^*))\in\nonegset$.

Note that, by construction, every string in $ \xy _\lambda$ contains an equal number of symbols $\wP _\lambda$ and $\wI _\lambda$.
Thus, $ \ltm _\lambda ( \xy _\lambda) =  \ltm ( \xy )$, since $\lambda$ and $\lambda^{-1}$ elements in $ \ltm _\lambda ( \xy _\lambda)$ cancel out; a consequence is that $\graph(\ltm_\lambda(\xy_\lambda))\in\nonegset$ is equivalent to $\graph(\ltm(\xy))\in\nonegset$.

Therefore, we only need to show that $\graph( \ltm _\lambda(\xx_\lambda \xy _\lambda^*))\in \nonegset$ and $\graph(\ltm _\lambda( \yy_\lambda\xy _\lambda^*))\in \nonegset$ if and only if $\lambda_{\min}\preceq \lambda \preceq \lambda_{\max}$.
Observe that
\[
\ltm_\lambda(\xx_\lambda\xy_\lambda^*)=\MC^*\lambda\MP\MC^*\ltm _\lambda ( \xy _\lambda)^* = \lambda\MC^*\MP\MC^* \ltm ( \xy)^*=\lambda\ltm(\xx\xy^*)
\]
and 
\[
\ltm_\lambda(\yy_\lambda\xy_\lambda^*)=\MC^*\lambda^{-1}\MI\MC^*\ltm _\lambda ( \xy _\lambda)^* = \lambda^{-1}\MC^*\MI\MC^* \ltm ( \xy)^* =\lambda^{-1}\ltm(\yy\xy^*).
\]
From~\Cref{pr:minimumcyclemean}, $\graph(\lambda \ltm (\xx \xy^* ))\in \nonegset$ iff 
\[
\lambda \preceq (\mbox{\normalfont mcm}( \ltm (\xx \xy^* )))^{-1},
\]
and $\graph(\lambda^{-1}\ltm (\yy \xy^* ))\in \nonegset$ iff 
\[
\lambda \succeq \mbox{\normalfont mcm}(\ltm (\yy \xy^* )) .
\qedhere\]
\end{proof}

We conclude this section with the proof of Theorem~\ref{th:correctness}.

\begin{proof}[Proof of Theorem~\ref{th:correctness}]
The proof is a consequence of~\Cref{th:extension4}.
Indeed, note that the "for" cycle of Algorithm~\ref{al:lambdasnn} computes \linebreak$ \ltm ( \xy ( \xx , \yy ,{\left \lfloor{\frac{n}{2}} \right \rfloor}))$, where $ \ltm ( \xx ) = \MC^* \MP \MC^*$, $ \ltm ( \yy ) = \MC^* \MI \MC^*$.
Therefore, $\lambda_{\min}$ and $\lambda_{\max}$, as defined in the previous theorem, coincide with $\lambda_{\min}$ and $\lambda_{\max}$ computed in Algorithm~\ref{al:lambdasnn}.
\end{proof}

\section{Conclusions}\label{se:conclusions}

In the present paper, we examine the \textit{Proportional-Inverse-Constant-Non-positive Circuit weight Problem} (PIC-NCP), which consists in finding all the values of $\lambda$ for which the parametric directed graph $\graphPIC$ does not contain circuits with positive weight.
The problem generalizes the NCP to a class of parametric graphs that is larger than the ones already studied in the literature.
After showing that the problem can be solved using linear programming, we present an algorithm that solves it in strongly polynomial time $\pazocal{O}(n^4)$ and provides a closed-formula expression for the lower and upper bound of the solution set.
The algorithm is based on a connection between square matrix operations in max-plus algebra, graph theory and formal languages.
The interest for this problem comes from the study of a specific class of discrete-event systems:
indeed, given a P-time event graph with at most $1$ initial token per place, the set of periods of all $d$-periodic trajectories that are consistent for the P-time event graph can be found, in strongly polynomial time, by solving a certain instance of the PIC-NCP%
~\cite{zorzenon2021periodic}.

To conclude, we state an open problem related to our work.
%
We remark that a more general class of NCP can be solved using the same algorithm presented in this paper.
Indeed, let us consider a max-plus Laurent polynomial in one variable $\lambda\in\R$ and matrix coefficients $A^{(-n_\MI)},A^{(-n_\MI+1)},\ldots,A^0,\ldots$ $A^{(n_\MP)}\in\Rmax^{n\times n}$, with $n_\MP,n_\MI\in\nato$:  $\bigoplus_{j=-n_\MI}^{n_\MP}\lambda^{j}A^{(j)}$.
The parametric precedence graph $\graph(\bigoplus_{j=-n_\MI}^{n_\MP}\lambda^{j}A^{(j)})$ can always be transformed into one of the form $\graphPIC$ by adding auxiliary nodes and arcs (obtained expanding $\lambda^j$ into $|j|$ products $\lambda\cdots\lambda$ if $j>1$ and $\lambda^{-1}\cdots\lambda^{-1}$ if $j<-1$).
However, the complexity for solving the NCP on $\graph(\bigoplus_{j=-n_\MI}^{n_\MP}\lambda^{j}A^{(j)})$
becomes pseudo-polynomial using this approach, since it increases with $n_\MP$ and $n_\MI$.
The problem can be solved in weakly polynomial time using linear programming; however, no strongly polynomial algorithm that solves it is known.
Its discovery would have \zorc{interesting} practical implications, as this algorithm could be used to check the existence of consistent $d$-periodic trajectories in strongly polynomial time complexity in general P-time event graphs (with no restriction on the number of initial tokens per place).

\section*{Acknowledgment}

This work was funded by the Deutsche Forschungsgemeinschaft (DFG, German Research Foundation), Projektnummer RA 516/14-1. Partially supported by the GACR grant 19-06175J, by  MSMT INTER-EXCELLENCE project\linebreak LTAUSA19098, by RVO 67985840, and by Deutsche Forschungsgemeinschaft (DFG, German Research Foundation) under Germany’s Excellence Strategy -- EXC 2002/1 "Science of Intelligence" -- project number 390523135 is also acknowledged.

\section*{Declaration of interest}
None.

\bibliographystyle{plainnat}
\bibliography{references}

\end{document}